\documentstyle[12pt]{article}
\title{Infinite Time Turing Machines}
\author{Joel David Hamkins{\footnotesize$^1$} and Andy Lewis}
\addtocounter{footnote}{1}\footnotetext{The research of the first author was
partially supported by a grant from the PSC-CUNY Research Foundation.}
 
\begin{document}


\newcommand{\wCK}{\omega_1^{\scriptscriptstyle CK}}
\newcommand{\ITleq}{\leq_{\scriptscriptstyle\infty}}
\newcommand{\ITlt}{<_{\scriptscriptstyle\infty}}
\newcommand{\ITequiv}{\equiv_{\scriptscriptstyle\infty}}
\newcommand{\Iff}{\mathbin{\longleftrightarrow}}
\newcommand{\converges}{\mathord{\downarrow}}
\newcommand{\iso}{\cong}
\newcommand{\zz}{Z$_0$}
\newcommand{\zzz}{Z$_0^{(0)}$}
\newcommand{\eth}{\hfill $\Box$\\}
\newcommand{\intersection}{\cap}
\newcommand{\sat}{\models}
\newcommand{\cross}{\times}
\newcommand{\prop}[1]{\noindent {\bf Proposition #1}}
\newcommand{\cor}[1]{\noindent {\bf Corollary #1}}
\newcommand{\lm}[1]{\noindent {\bf Lemma #1}}
\newcommand{\class}{\noindent {\bf Class Participation}: }
\newcommand{\assign}[1]{\noindent{\bf Assignment E.#1}}
\newcommand{\th}[1]{\noindent {\bf Theorem #1}}
\newcommand{\df}[1]{\noindent {\bf Definition #1}}
\newcommand{\nh}{\newline \hspace*{.5in}}
\newcommand{\un}[1]{$\underline{#1}$}
\newcommand{\uA}{$\underline{A}$}
\newcommand{\uR}{$\underline{R}$}
\newcommand{\uN}{$\underline{N}$}
\newcommand{\uZ}{$\underline{Z}$}
\newcommand{\uQ}{$\underline{Q}$}
\newcommand{\be}{\begin{enumerate}}
\newcommand{\ee}{\end{enumerate}}
\newcommand{\bd}{\begin{description}}
\newcommand{\ed}{\end{description}}
\newcommand{\lub}{$ lub $}
\newcommand{\ran}{$ ran $}
\newcommand{\card}[1]{$\overline{\overline{#1}}$}
\newcommand{\ls}{\mbox{$<$}}
\newcommand{\rs}{\mbox{$>$}}
\newcommand{\rl}{\mbox{$\omega^{\omega}$}}
\newcommand{\oi}{\mbox{is order isomorphic to }}
\newcommand{\TS}{\mbox{topological space }}
\newcommand{\can}{\mbox{$2^{\omega}$}}
\newcommand{\alf}{\mbox{$\alpha$}}
\newcommand{\seq}{\mbox{$\omega^{< \omega}$}}
\newcommand{\canseq}{\mbox{$^{< \omega}2$}}
\newcommand{\calP}{\mbox{${\cal P}$}}
\newcommand{\calU}{\mbox{${\cal U}$}}
\newcommand{\calA}{\mbox{${\cal A}$}}
\newcommand{\calB}{\mbox{${\cal B}$}}
\newcommand{\calT}{\mbox{${\cal T}$}}
\newcommand{\calS}{\mbox{${\cal S}$}}
\newcommand{\Lx}{\mbox{${\bf L_\omega_1^x}$}}
\newcommand{\kn}{\mbox{$[\kappa]^{n}$}}
\newcommand{\kw}{\mbox{$[\kappa]^{<\omega}$}}
\newcommand{\ka}{\mbox{$[\kappa]^{\alpha}$}}
\newcommand{\T}{\mbox{Trace}}
\newcommand{\B}{\mbox{$\Pi^1_2\,\,\,$}}
\newcommand{\A}{\mbox{$\Pi^1_1\,\,\,$}}
\newcommand{\Tpro}{\mbox{$p[T]$}}
\newcommand{\Spro}{\mbox{p[S]$}}
\newcommand{\TTpro}{\mbox{$p[T_2]$}}
\newcommand{\hlms}{\mbox{\rule{.9em}{.9em}}}
\newcommand{\ot}{\mathop{\hbox{\rm ot}}}
\newcommand{\st}{\mid}
\def\underTilde#1{{\baselineskip=0pt\vtop{\hbox{$#1$}\hbox{$\sim$}}}{}}
\def\undertilde#1{{\baselineskip=0pt\vtop
  {\hbox{$#1$}\hbox{$\scriptscriptstyle\sim$}}}{}}
\font\am=msam10 at 12pt
\font\bm=msbm10 at 12pt
\newcommand{\restrict}{\mathbin{\hbox{\am\char'26}}}
\def\concat{\mathbin{{\,\hat{ }\,}}}
\def\force{\mathbin{\hbox{\am\char'15}}}
\def\notforce{\mathbin{\hbox{\bm\char'61}}}
\def\TC{\mathop{\hbox{\sc tc}}\nolimits}
\def\R{\hbox{\bm R} }
\def\Z{\hbox{\bm Z} }
\def\Q{\hbox{\bm Q} }
\def\N{\hbox{\bm N} }
\def\I{\hbox{\bm I} }
\def\ZN{\hbox{\bm Z$^{>0}$ }}
\def\xv{\hbox{\bf x} }
\def\yv{\hbox{\bf y} }
\def\zv{\hbox{\bf z} }
\def\df{\bf}
\def\lt{\mathrel{\mathchoice{\scriptstyle<}{\scriptstyle<}
   {\scriptscriptstyle<}{\scriptscriptstyle<}}}
\def\from{\mathbin{\vbox{\baselineskip=3pt\lineskiplimit=0pt
                         \hbox{.}\hbox{.}\hbox{.}}}}
\def\tlt{\triangleleft}
\def\theorem{\claimtype=0\sayclaim}
\def\<#1>{\left\langle#1\right\rangle}
\def\set#1{\{\,#1\,\}}
\newcommand{\union}{\cup}
\font\nineam=msam9
\newcommand{\jump}{{\!\hbox{\nineam\char'117}}}
\newcommand{\Jump}{{\!\hbox{\nineam\char'110}}}
\newcommand{\proof}{\noindent{\bf Proof: }}

\maketitle
 
$$\vbox{\hsize=5in \setlength{\baselineskip}{1pt} 
{\small 
\noindent {\bf Abstract.}
We extend in a natural way the operation of Turing machines to infinite
ordinal time, and investigate the resulting supertask theory of
computability and decidability on the reals. Every $\Pi^1_1$ set, for
example, is decidable by such machines, and the semi-decidable sets
form a portion of the $\Delta^1_2$ sets. Our oracle concept leads to a
notion of relative computability for sets of reals and a rich degree
structure, stratified by two natural jump operators.}}$$

\noindent In these days of super-fast computers whose speed seems to be
increasing without bound, the more philosophical among us are perhaps
pushed to wonder: what could we compute with an {\it infinitely} fast
computer? By proposing a natural model for supertasks---computations
with infinitely many steps---we provide in this paper a theoretical
foundation on which to answer this question. Our model is simple: we
simply extend the Turing machine concept into transfinite ordinal time.
The resulting machines can perform infinitely many steps of
computation, and go on to more computation after that.  But
mechanically they work just like Turing machines.  In particular, they
have the usual Turing machine hardware; there is still the same smooth
infinite paper tape and the same mechanical head moving back and forth
according to a finite algorithm, with finitely many states.  What is
new is the definition of the behavior of the machine at limit ordinal
times. The resulting computability theory leads to a notion of
computation on the reals, concepts of decidability and
semi-decidability for sets of reals as well as individual reals, two
kinds of jump-operator, and a notion of relative computability using
oracles which gives a rich degree structure on both the collection of
reals and the collection of sets of reals. But much remains unknown; we
hope to stir interest in these ideas, which have been a joy for us to
think about.

\newpage
 
There has been much work in higher recursion theory analyzing
well-founded computations on infinite objects (see e.g. \cite{Sa}).
Much of that theory, however, grows out of the analogy which associates
the $\Delta^1_1$ sets with the finite sets and the $\Pi^1_1$ sets with
the semi-decidable sets. It therefore gives
a different analysis than ours, since in our account, the $\Pi^1_1$ sets will
become actually {\it decidable}, along with their complements the 
$\Sigma^1_1$ sets, and the semi-decidable sets, with the jump operator, 
will stratify the class of $\Delta^1_2$ sets.

Various philosophers and physicists have investigated supertasks, 
or tasks involving infinitely many steps, the first of which is done,
for example, in a half second, the next in a quarter second, and so on,
so that the entire job is complete in a finite amount of time.
Thomson's lamp \cite{Thom}, for example, is on between time $t=0$ and
$t=1/2$, off until $t=3/4$, on until $t=7/8$, and so on. More useful
supertasks, perhaps, have been proposed which determine the truth of an
existential number-theoretic question, such as whether there are
additional Fermat primes, by ever more rapidly checking the
instances of it so that they are all checked in a finite amount of
time. What is intriguing about the physicist's analysis of supertasks
is that they have been able to construct general relativistic models in
which the supertasks can apparently be carried out \cite{Ear1},
\cite{Ear2}, \cite{Pit}, \cite{Hog92}, \cite{Hog94}. The models generally
involve an agreement between two observers, one of whom performs the
rote steps of computation looking for a counterexample, the other of
whom, while flying ever more rapidly around the first observer, waits
patiently for what to him will be a finite amount of time for a signal
that a counterexample has been found. Earman shows how it can happen
that the entire infinite past half-life of one observer is in the
causal past of a point in the life of another observer. Hogarth
\cite{Hog94} discusses more complicated arrangements in which the truth
of any arithmetic statement can be computed. But as we will show in
this paper, the supertask concept allows one to compute the truth of
even more complicated statements than this. 
What we are interested in here is not so much finding what is 
physically possible to compute in a supertask so much as what is 
{\it mathematically} possible. Though the physicists may explain 
how it is possible to carry out a supertask in a finite amount of 
time, we, being focused on the algorithm, will nevertheless regard
the supertask computations as being infinite in the sense that they 
involve infinitely many steps of computation. 

A word of credit is due. Jeffrey Kidder defined infinite time Turing
machines in 1989, and he and the first author of this paper worked out
the early theory while they were graduate students together at the
University of California in Berkeley. The notes sat neglected,
unfortunately, in the filing cabinet of the first author, carted from
one university to another, but were never forgotten. With a fresh look in
1996, a few of the formerly puzzling questions were answered, and a
more advanced theory developed, which we give in this paper.
 
We will assume complete familiarity with the notions of Turing machines
and ordinals and, in describing our algorithms, take the high road to
avoid getting bogged down in Turing machine minutiae.  We hope the readers
will appreciate our saving them from reading what would otherwise
resemble computer code.

We begin by describing in section 1 how the new machines work. Next, 
in section 2, we investigate their power, which lies between 
$\Pi^1_1$ and $\Delta^1_2$, and generalize the classical s-m-n and 
Recursion theorems to the supertask context. Our analysis of the 
lengths of supertask computations leads in section 3 to the 
notions of clockable and writable ordinals. We present in section 4 
the supertask halting problems, a lightface and boldface version, and 
prove the Lost Melody Theorem, which has the intriguing consequence 
that there are noncomputable functions whose graphs are infinite time 
decidable. This leads naturally in section 5 to the notion of oracles 
and the two jump operators corresponding to the two halting problems. 
We develop in section 6 the basic features of the structure of 
infinite time degrees; a key feature of this structure, having no 
classical analog, is the interplay between reals as oracles and 
sets of reals as oracles. We highlight, in section 7, the descriptive 
set-theoretic aspects of our theory and, in section 8, the 
connections with admissible set theory; it turns out that the 
supremum of the writable ordinals is a rather large countable 
ordinal, being recursively inaccessible and indeed indescribable by 
semi-decidable properties. We conclude the paper by giving five 
equivalences to the question of whether every clockable ordinal is 
writable, a question which remains open.
 
\section{How the Machines Work}
 
Let us now describe exactly how an infinite time Turing machine works.
Like a Turing machine, our supertask machines have a head which moves
mechanically back and forth, reading and writing $0$s and $1$s on a
tape according to a finite algorithm $p$. For convenience, we will set
up our machines with three separate tapes, one for input, one for
scratch work, and one for output.  
 
\vspace*{3in}
 
\special{eps: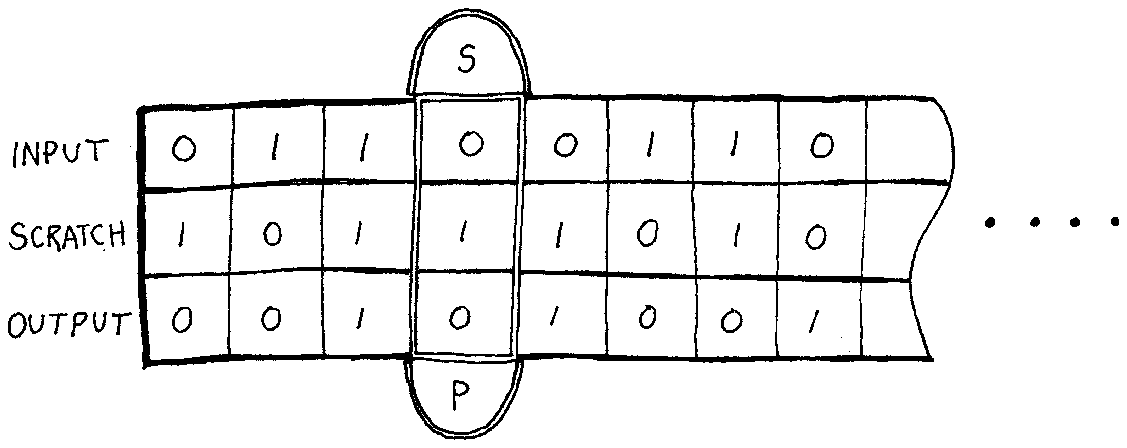 x=5.6in. y=2.8in.}
 
\noindent The machine begins, like a
Turing machine, with the head resting in anticipation on the first cell
in a special state called the {\it start} state. The input is written
on the input tape, and the scratch tape and the output tape are filled
with zeros. At each step of computation, the head reads the cell values
which it overlies, reflects on its state, consults the program about
what should be done in such a situation and then carries out the
instructions: it writes a $0$ or $1$ on (any of) the tapes, moves left or right
and switches to a new state accordingly. This procedure determines the
configuration of the machine at stage $\alpha+1$, given the
configuration at stage $\alpha$, for any $\alpha$. It remains to
somehow take a limit of these computations in order to identify the
configuration of the machine at stage $\omega$ and, more generally, at
limit ordinal stages in the supertask computation. 
To set up such a limit ordinal configuration, the
head is plucked from wherever it might have been racing towards, and
placed on top of the first cell. And it is placed in a special
distinguished {\it limit} state.  Now we need to take a limit of the
cell values on the tape. And we will do this cell by cell according to
the following rule:  if the values appearing in a cell have converged,
that is, if they are either eventually $0$ or eventually $1$ before the
limit stage, then the cell retains the limiting value at the limit
stage. Otherwise, in the case that the cell values have alternated from
$0$ to $1$ and back again unboundedly often, we make the limit cell
value $1$. This is equivalent to making the limit cell value the
$\limsup$ of the cell values before the limit.  This completely
describes the configuration of the machine at any limit ordinal stage
$\beta$, and the machine can go on computing to $\beta+1$, $\beta+2$,
and so on, eventually taking another limit at $\beta+\omega$ and so on
through the ordinals. If at any stage, the machine finds itself in the
special {\it halt} state, then computation ceases, and whatever is
written on the output tape becomes the official output.  Otherwise, the
infinite time machine will compute endlessly as the ordinals fall one
after another through the transfinite hourglass.
 
In this way every infinite time Turing machine program $p$ determines a
function. On input $x$, we can run the machine with program $p$, and,
if it halts, there will be some output, which we denote by
$\varphi_p(x)$. The domain of $\varphi_p$ is simply the collection of
$x$ which lead to a halting computation. Notice that the natural input
for these machines is an {\it infinite} binary string $x\in 2^\omega$.
Thus, the infinite time computable functions are partial functions on
Cantor space. In this paper, we will refer to the elements of Cantor
space $2^\omega$ as `reals', and think of the computable functions as
functions on the reals. In particular, we will denote $2^\omega$ by
$\R$.  By adding extra input tapes, we may have functions of more than
one argument, and, in the usual Turing machine argument, these can be
simulated by, for example, interleaving the digits of the inputs and
using a machine with only one input tape. Also, we will regard $0$ and
$1$ as elements of $\R$, by the convention, for example, in which $0$
represents $\<0,0,0,\ldots>$ and $1$ represents $\<1,0,0,0,\ldots>$.
 
Let us now make a few basic definitions. A partial function $f\from \R^k
\rightarrow \R$ is infinite time {\it computable} when there is a
program $p$ such that $f=\varphi_p$. For simplicity, we will assume by
some suitable coding mechanism that a program is represented by a
natural number, and furthermore, that every natural number represents a
program. A set of reals $A$ is infinite time {\it decidable} when the
characteristic function of $A$ is infinite time computable (and for the sake
of brevity, we will normally use just the term ``decidable," unless we
fear some misunderstanding).  The set $A$ is infinite time
{\it semi-decidable} when the function which gives the affirmative values,
the function with domain $A$ and constant value $1$, is infinite time
computable. Thus, a set is semi-decidable exactly when it is the domain
of a computable function, since it is a simple matter to modify a
program to change the output to the constant $1$.  We can also
stratify the computable sets according to how long the computations
take. Thus, a set is $\alpha$-decidable when the characteristic
function of the set is computable by a machine which on any input takes
fewer than $\alpha$ many steps. Thus, restricting to the case of finite
input and time, a function $f\from 2^{\lt \omega} \rightarrow 2^{\lt
\omega}$ is $\omega$-computable exactly when it is computable in the
Turing machine sense.
 
Now that we have made the initial definitions, let us get started by
proving that we may restrict our attention always to the countable
ordinals.
 
\newtheorem{countable}{Theorem}[section]
\begin{countable} 
Every halting infinite time computation is countable. 
\label{onetwo}
\end{countable}
 
\proof  Suppose the supertask computation of program $p$ on input $x$
has run for uncountably many steps without halting. We will show that
this computation will never halt.  Let us say that a {\it snapshot} of
a computation is a complete description of the configuration: it
specifies the program being used, the state and position of the head
and the complete contents of each of the tapes. All this information
can be coded in some canonical way into a real.  Consider the snapshot
of the $\omega_1$-stage of computation. We will argue that in fact this
very same snapshot occurred earlier, at some countable stage, and that
the computation is caught in an infinite loop which will repeat
forever. First, observe that at stage $\omega_1$ the head is on the
first cell in the limit state, as it is at any limit ordinal stage.
Next, observe that if the value of any cell is $0$ at stage $\omega_1$,
then there must be some countable stage at which the cell had the value
$0$ and never subsequently changed. If the cell has the value $1$ at
stage $\omega_1$, then there are two possibilities: either at some
stage the cell obtained a value of $1$ and was not subsequently changed
to a $0$, or else the value of the cell alternated unboundedly often
before stage $\omega_1$. Now we apply a simple cofinality argument.
Since there are only countably many cells, by taking a countable supremum we
can find a stage $\alpha_0$ where the cells which eventually
stabilize have already all stabilized. After this stage, the only cells
which change are the ones which will change cofinally often, and so
there must be a sequence of countable ordinals
$\alpha_0<\alpha_1<\alpha_2<\cdots$ such that between $\alpha_n$ and
$\alpha_{n+1}$ all the cells which change at all after $\alpha_n$ have
changed value at least once by $\alpha_{n+1}$. Let $\delta = \sup_n
\alpha_n$. At this limit stage, the head is on the first cell in the
limit state, the cells which stabilize before $\omega_1$ have
stabilized before $\delta$, and the cells which change values
unboundedly often before $\omega_1$ have changed values unboundedly
often before $\delta$. Thus, the snapshot at stage $\delta$ is the same
as the snapshot at stage $\omega_1$. Thus, the computation has repeated
itself. And moreover it has done so in a very strong way: none of the cells
which are zero at stage $\delta$ will ever again turn to one, so the 
computation will cycle endlessly. Limit stages in which the repeating snapshot
has repeated unboundedly often will again 
be the very same repeating snapshot, and so the computation is caught
in an infinite repeating loop.\eth
 
Let us point out that it is possible for a computation to repeat---for
the very same snapshot to occur twice during a computation---but for
the computation nevertheless eventually to escape this infinite loop.
This could occur, for example, if the limit of the repeating snapshots
is not the same snapshot again. Thus, after repeating $\omega$ many
times, the program could escape the loop and go on to compute something
else or even halt.  Let us say, then, officially, that a computation
repeats itself only when the very same snapshot occurs at two limit
ordinal stages, and that between these stages the cells which are $0$
at the limit never subsequently turn to $1$ (we allow the $1$s to turn
to $0$ and back again).  This is equivalent to requiring that the limit
of the repeating snapshots is again the very same repeating snapshot.
Such computations, therefore, are truly caught in an infinite loop.
 
\newtheorem{cycle}[countable]{Corollary}
\begin{cycle}
Every infinite time computation either halts or repeats itself in
countably many steps.
\end{cycle}
 
\proof The previous proof shows that if a computation does not halt,
then the snapshot of the machine configuration at stage $\omega_1$
appears earlier, and, in fact, appears unboundedly often before
$\omega_1$, since $\alpha_0$ can be chosen arbitrarily large below
$\omega_1$. Moreover, since these snapshots occurred beyond $\alpha_0$,
none of the $0$s in the repeating snapshot ever again turns to $1$. So
a non-halting computation repeats itself by some countable stage.\eth
 
\section{The Power of Infinite Time Machines}
 
How powerful are these machines? Perhaps the first thing to notice is
that the halting problem for Turing machines is infinite time
decidable. This is true because with an infinite time Turing machine
one can simulate an ordinary Turing machine computation. Either the
simulation halts in finitely many steps, or else after $\omega$ many
steps the machine reaches the limit state, and so by giving the output
Yes or No, respectively, in these two situations, the halting problem is
solved. Thus infinite time Turing machines are more powerful than
ordinary Turing machines:  they can decide sets which are undecidable
by Turing machines.  The next theorem greatly improves on this.
 
\newtheorem{arithtruth}{Theorem}[section]
\begin{arithtruth}
The truth of any arithmetic statement is infinite time decidable.
\label{arith}
\end{arithtruth}
 
\proof It is easy to see by an inductive argument on formulas
that first order arithmetic truth is decidable: given a
statement $\exists n\varphi(n,\vec x)$, one simply directs the machine
to try out all the possible values of $n$ and test the truth of
$\varphi(n,\vec x)$.\eth
 
But we can do much better even than this. In the next theorem we will
introduce the first argument which really seems to use the full
computational power of these machines. A relation $\triangleleft$ on a subset
of $\omega$ can be coded by the real $x$ such that $x(\<n,k>)=1$ exactly
when $n\triangleleft k$, where $\<\cdot,\cdot>$ is some canonical pairing
function. In this way every real $x$ codes some relation
$\triangleleft$. Let $WO$ be the set of reals coding well-orders.
It is well known that WO is a complete $\Pi^1_1$ set, 
in the sense that if $A$ is another $\Pi^1_1$ set then there is a 
recursive function $f$ on reals such that $x\in A\iff f(x)\in WO$.
 
\newtheorem{decwo}[arithtruth]{Count-Through Theorem}
\begin{decwo}
$WO$ is infinite time decidable.
\label{twotwo}
\end{decwo}
 
\proof  Suppose we are given a real $x$ on the input tape of 
an infinite time Turing machine.  We will describe a supertask 
algorithm which will determine if this
real codes a well order. Certainly in $\omega$ many steps we can
determine if the real codes a relation which is reflexive, transitive,
and antisymmetric.  This involves simply checking that if $x$ says that
$n$ is related to $k$ and that $k$ is related to $r$, then it also says
that $n$ is related to $r$, and so on. Every instance of this can be
systematically checked.  Thus, we may assume that the real $x$ survives 
this first test, and therefore codes a
linear order, before continuing further. Next, in $\omega$ many steps,
we can find the particular natural number $n$ which is the least
element in the relation coded by $x$.  This can be done by first, 
keeping a current guess written on the scratch tape, updating it every
time a number is found which precedes it in the relation coded by $x$,
and second, flashing a flag on and then off again every time we change
the guess. In the limit, if the flag is on, it means that we changed
our guess infinitely often, and thus the real $x$ does not code a
well-order.  If the flag is not on, then it must be that the minimal
element $n$ of the relation is sitting on the scratch tape, having
survived all challengers. In another $\omega$ many steps, we can 
go through the relation and erase all mention of the number $n$ from
the field of the relation. This produces a real coding a relation with
a smaller field. Now we simply iterate this. That is, in $\omega$ many
steps, first find the least element of the relation coded by the real
currently on the tape, and in another $\omega$ many steps, erase all
mention of this element, and repeat with the new smaller relation.
There is a slight complication at the compound limits (limits of
limits) since at such stages we will have a lot of garbage on the
scratch tape, but because we were gradually erasing elements from the
field of the relation coded by the real on the input tape, the input
tape has stabilized to the intersection of those relations, which is
exactly what we want there. By flashing a flag on and then off again
every time we reach a limit stage, we can recognize a compound limit as
a limit stage in which this flag is on, and then in $\omega$ many steps
wipe the scratch tape clean before continuing with the algorithm.  If
the original real codes a relation which is not a well-order, then
after its well-founded part has been erased, there will come a stage
when it has no least member, and the machine will discover this, since
the guesses for the least member at that stage will not converge. If
the real does code a well-order, then the elements of the field will
gradually be erased until the machine finds that the field of the
relation is empty. Thus, in any case the machine will know whether the
real codes a well-order, and so $WO$ is decidable.\eth
 
\newtheorem{decpi}[arithtruth]{Corollary}
\begin{decpi}
Every $\Pi^1_1$ set is infinite time decidable. Hence, every
$\Sigma^1_1$ set is infinite time decidable.
\end{decpi}
 
\proof  It is well-known that every $\Pi^1_1$ set $A$ reduces to $WO$
in the sense that there is some recursive function on reals $f$ such
that $x\in A\leftrightarrow f(x)\in WO$. Furthermore, it is clear that
an infinite time machine can compute any recursive function on reals.
So, fix $A$ and $f$, and consider the algorithm which on input $x$
first computes $f(x)$, and then determines if $f(x)\in WO$. This machine
decides whether $x$ is in $A$.\eth
 
The collection of decidable sets extends further up the analytical
hierarchy. A set $A$ is said to be $\beta{-}\Pi^1_1$, where $\beta$ is
a recursive ordinal coded by some recursive relation $E$ on $\omega$,
when for each $\alpha\leq\beta$ there is a set $A_\alpha$, with
$A_\beta=\emptyset$, such that, first, $\set{(k,x)\mid x\in A_{|k|}}$
is $\Pi^1_1$, where $|k|$ denotes the order-type of $k$ with respect to
$E$, and, second, that $x\in A$ iff there is an odd $\alpha<\beta$ such
that $x\in\intersection_{\delta<\alpha}A_\delta\setminus A_\alpha$ (see
\cite{D}). This class of sets extends beyond the $\Pi^1_1$ sets.  We
can extend it still further by allowing more complicated relations
$E$.  Specifically, let us say that a real is {\it writable} when there is an
infinite time Turing machine which can write it as the final output on
input $0$.  An ordinal will be regarded as writable when there is a
writable real coding that ordinal. We will show later that such
ordinals extend far beyond the recursive ordinals.  Indeed, their
supremum is recursively inaccessible.  We can naturally extend the
definition of the $\beta{-}\Pi^1_1$ sets to the situation when $\beta$
is a writable ordinal.

Suppose now that $A$ is $\beta{-}\Pi^1_1$ where $\beta$ is a writable
ordinal.  Consider the algorithm which first writes the relation $E$ coding
$\beta$ on a portion of the scratch tape, and then, using the $\Pi^1_1$
algorithm makes a list of which $n$ have the property that the input
$x$ is in $A_{|n|}$. Finally, by counting through the relation coding
$\beta$, the algorithm searches for an odd ordinal $\alpha<\beta$ such
that $x\in\intersection_{\delta<\alpha}A_\delta\setminus A_\alpha$. This
algorithm will decide whether $x$ is in $A$, giving the following
corollary.
 
\newtheorem{decbeta}[arithtruth]{Corollary}
\begin{decbeta}
If $\beta$ is a writable ordinal, then every $\beta{-}\Pi^1_1$ set is decidable.
\end{decbeta}
 
We aim now to establish a limit on the complexity of decidable sets.
In the previous section, we introduced the idea of a snapshot,
which is a real which codes the complete description of an infinite time
Turing machine while
it is computing. Thus, in some canonical manner, a snapshot codes the
program the machine is running, the position and state of the head and
the complete contents of each of the three tapes. Now let us say that a
transfinite sequence of snapshots {\it accords} with the program $p$
when each successive snapshot is obtained by running the program $p$ on
the configuration described by the previous snapshot and the limit
snapshots are obtained from the earlier ones according to the
computation rules. We will say that a sequence of snapshots according
to a program is {\it settled} when the last snapshot has either
obtained a halting state or else repeats an earlier snapshot, in the strong
sense of a computation repeating itself used after Theorem \ref{onetwo}.  
The sequence of snapshots represents in the first case a halting
computation and in the second a computation which will endlessly
repeat. Thus, a settled sequence of snapshots informs us of the outcome
of an infinite time computation. With these ideas we can establish the
complexity of such computations.
 
\newtheorem{graphdlta}[arithtruth]{Complexity Theorem}
\begin{graphdlta}
The graph of every infinite time computable function is $\Delta^1_2$.
Hence, every decidable set and, indeed, every semi-decidable set, is
$\Delta^1_2$.
\end{graphdlta}
 
\proof Suppose that $f$ is an infinite time computable function, computed by
the program $p$. Now observe that $f(x)=y$ if and only if there is a
real $z$ coding a well-ordered sequence of snapshots according to $p$
on input $x$ with output $y$. Thus, the
graph of $f$ is $\Sigma^1_2$.  For the other half, remember that by
Theorem \ref{onetwo} every supertask computation either halts or repeats
in countably many steps, and that therefore we only need to consider
the settled sequences of snapshots.  That is, $f(x)=y$ if and only if
every real $z$ coding a well-ordered settled sequence of snapshots
according to $p$ on input $x$ shows the output as $y$. Thus the graph
of $f$ is also $\Pi^1_2$, and therefore $\Delta^1_2$, as desired.\eth
 
It follows, of course, that every co-semi-decidable set is also
$\Delta^1_2$. Later we will show, using the halting problem for infinite
time computations, that there are semi-decidable sets which are not 
decidable, and therefore also co-semi-decidable sets which are
not semi-decidable.  Therefore, the semi-decidable sets form a proper
subclass of the $\Delta^1_2$ sets, as illustrated in the diagram
below.

\vspace*{4.5in}
 
\special{eps: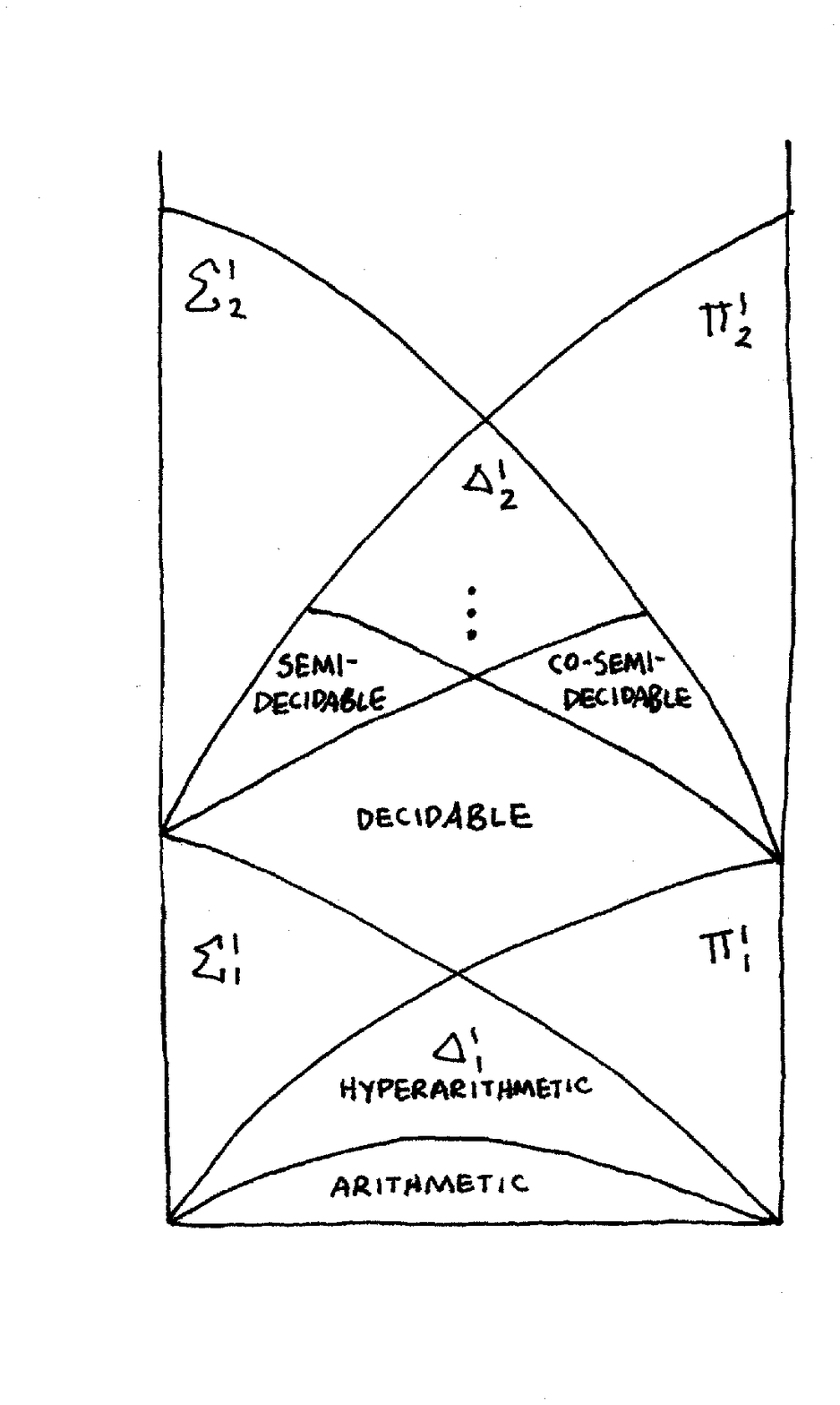 x=3.6in. y=4.5in.}
 
\noindent In the later sections we will show how the jump operators
actually stratify the class of $\Delta^1_2$ sets. For now, let us 
exactly identify the classes of sets which are decidable by algorithms which
need to take relatively few limits. 
 
\newtheorem{arithfin}[arithtruth]{Theorem}
\begin{arithfin}
The arithmetic sets are exactly the sets which can be decided by an
algorithm using a bounded finite number of limits.
\end{arithfin}
 
\proof The inductive argument in Theorem \ref{arith} essentially shows
that a $\Sigma^0_n$ set can be decided by an algorithm using at most $n$
limits. For the converse direction, let us define the notation 
$\varphi_{p,\alpha}(x)$ to mean the snapshot of the computation of 
$\varphi_p(x)$ after exactly $\alpha$ many steps of computation. 
We claim that the relation ``s is an initial segment of 
$\varphi_{p,\omega n}(x)$,'' that is, $s\subset\varphi_{p,\omega n}(x)$, 
is arithmetic. Certainly this is true 
when $n=0$. Now consider $n+1$. Let $p'$ be the program which 
does the same thing when starting that $p$ does at a limit stage (thus, 
for example, $p'$ could simply at the first step change to the 
{\it limit} state, and thereafter mimic $p$). The only way there
can be a $0$ on a cell of the tape in the limit $\omega(n+1)$ is when the 
configuration is set up at stage $\omega n$ so that the cell in question 
will be $0$ from some point on. Thus, 
$s\subset\varphi_{p,\omega(n+1)}(x)$ if and only if for every $i<\ell(s)$
which corresponds to the value in a cell, 
$s(i)=0\Iff\exists N\,\forall m{>}N\,\forall t\, [t\subset\varphi_{p,\omega n}(x)\wedge
\ell(t)>m\longrightarrow \varphi_{p',m}(t)(i)=0]$. The very last clause is 
equivalent to $\varphi_{p,\omega n+m}(x)(i)=0$, because of the definition
of $p'$ and the assumption on $t$. Thus, the relation is arithmetic at 
$\omega(n+1)$, and the claim is proved. The theorem follows, since if 
$A$ is decided always in $n$ limits by program $p$, then $x\in A$ if and 
only if $\<1>\subset\varphi_{p,\omega n}(x)$, and so $A$ is arithmetic.\eth
 
\newtheorem{hyperarith}[arithtruth]{Theorem}
\begin{hyperarith}
The hyperarithmetic sets, the $\Delta^1_1$ sets, are exactly the sets which 
can be decided in some bounded recursive ordinal length of time.
\end{hyperarith}
 
\proof Suppose that $A$ is decidable by program $p$, and the
computation $\varphi_p(x)$ always halts before the recursive ordinal
$\alpha$. We want to prove that $A$ is $\Delta^1_1$. The ordinal
$\alpha$ is coded by some recursive relation $\triangleleft$ on
$\omega$, coded by the real $y$.  Thus, $x\in A$ if and only if there
is a sequence of snapshots of length at most $\alpha$ according to $p$
on input $x$ which shows the computation to halt with output $1$. This
is a $\Sigma^1_1$ property, since whether a real codes a sequence of
length at most $\alpha$ is equivalent to whether there is an order
isomorphism from the given relation to an initial segment of the
relation coded by $y$, and $y$ itself is recursive. Thus, $A$ is
$\Sigma^1_1$.  We also know, however, that $x\in A$ if and only if
every real which codes a sequence of snapshots according to $p$ on
input $x$ which is well-founded at least to $\alpha$, shows the
computation to halt with output $1$. This shows that $A$ is $\Pi^1_1$,
and so $A$ is $\Delta^1_1$, as desired.
 
Conversely, suppose that $A$ is $\Delta^1_1$. We want to show that $A$
is infinite time decidable by an algorithm taking some bounded recursive
length of time. Every hyperarithmetic set has a recursive Borel code, a
recursive well-founded tree where each node is labeled with instructions for
building the set from the earlier nodes by taking unions,
intersections, or complements, and the minimal nodes are labeled with basic
open sets (the whole set $A$ being represented by the top 
node).  Consider the algorithm which, on input $x$, systematically works
through this tree and keeps track of whether $x$ is in or out of the
set coded by each particular node. Thus, for example, the algorithm
determines whether $x$ is in each of the open sets at the bottom 
of the
tree, and if a node is labeled with intersection, then the algorithm
determines whether $x$ is in the set represented by that particular
node by checking whether it has already said that $x$ is in each of the
sets represented by the earlier nodes. The top node represents $A$, and
so this algorithm will decide whether $x$ is in $A$ or not.
Furthermore, since the tree is recursive, the height of the tree is
recursive, and so this algorithm will take a bounded recursive length
of time to determine the final answer. If the nodes are handled according 
to their rank in the tree, since it takes $\omega$ many steps to handle
each node, the whole algorithm will take $\rho\omega$ many steps, where 
$\rho$ is the rank of the tree. Since the tree is recursive, this is 
a recursive ordinal.\eth
 
Let us conclude this section by proving that two results from classical
recursion theory hold also for infinite time Turing machines.
 
\newtheorem{smn}[arithtruth]{The s-m-n Theorem}
\begin{smn}
There is a primitive recursive function $s$ defined on the natural numbers 
such that $\varphi_p(\vec k,\vec x)=\varphi_{s(p,\vec k)}(\vec x)$.
\end{smn}
 
\proof  The classical proof works in this context just as well. One needs
only to check that there is a simple uniform procedure to convert a 
program $p$ which computes a function of arity $m+n$ into a program 
$s(p,\vec k)$ which computes the same function with the first $m$ arguments
fixed as $\vec k$, thereby computing a function of arity $n$.\eth
 
\newtheorem{recthm}[arithtruth]{The Recursion Theorem}
\begin{recthm}
For any infinite time computable total function $f:\N\to\N$, 
there is a program $p$ such that $\varphi_p=\varphi_{f(p)}$.
\label{Recthm}
\end{recthm}
 
\proof  The classical proof works also for this theorem. Fix the 
function $f$. Let $r$ be the program which computes the function 
$\varphi_r(q,x)=\varphi_{f(s(q,q))}(x)$. Let $p=s(r,r)$, and observe 
that $\varphi_p(x)=\varphi_{s(r,r)}(x)=\varphi_r(r,x)=\varphi_{f(s(r,r))}(x)=
\varphi_{f(p)}(x)$, as desired.\eth
 
So that the reader does not think that all the classical arguments will
easily generalize, let us just mention that in the supertask 
context there are noncomputable functions whose graphs are
semi-decidable. Indeed, there is a total constant function which is not
computable, but whose graph is nevertheless decidable! This will be
proved at the end of section four.
 
\section{Clockable Ordinals}
 
The study of infinite time computations leads inexorably to the desire
to know exactly how long the computations take. And this leads to the
notion of clockable ordinals. An ordinal $\alpha$ is {\it clockable}
when there is a program which on input $0$ halts in exactly $\alpha$
many steps of computation. That is, the $\alpha^{\rm th}$ step of
computation is the act of changing to the {\it halt} state.  Any
natural number $n$, for example, is clockable---one simply uses a
machine which counts through $n$ states before halting. The ordinal
$\omega$ is also clockable, since a machine can be programmed, for
example, to move the head always to the right, until a limit state is
obtained, and then halt. It is not difficult to prove that if $\alpha$
is clockable, then so are $\alpha+1$ and $\alpha+\omega$. These simple
ideas show that every ordinal up to $\omega^2$ is clockable, and
moreover, if $\alpha$ is clockable, then so is $\alpha+\beta$ for any
$\beta<\omega^2$.  Let us now complete the warm-up for this section by
proving that $\omega^2$ itself is also clockable. Since $\omega^2$ is
the first ordinal which is a limit of limit ordinals, we will simply
design a machine which can recognize such compound limits.  At each
limit stage, let the machine flash a flag on and then off again, and
then simply wait for the next limit stage. The stage $\omega^2$ will be
the first time that the machine is in a limit state and sees that the
flag is on.  (of course, one needs to put this flag, signaling that we
are done, on the very first cell of the tape, in order that the machine
will not need extra steps after $\omega^2$).  Direct the machine to halt
when this occurs. Thus, $\omega^2$ is clockable.  The reader may enjoy
writing programs to halt at $\omega^3$ or
$\omega^{\omega^2+5}$.  The next theorem shows that clockable ordinals
are plentiful.
 
\newtheorem{recclock}{Recursive Clocks Theorem}[section]
\begin{recclock}
Every recursive ordinal is clockable.
\label{threeone}
\end{recclock}
 
\proof Suppose that $\alpha$ is a recursive ordinal which we
would like to show is clockable. By the remarks of the previous
paragraph, we may assume that $\alpha$ is a compound limit ordinal,
since any ordinal above $\omega^2$ is the sum of a compound limit
ordinal with an ordinal below $\omega^2$.  We will describe a supertask
algorithm which will halt in exactly $\alpha$ many steps. We will do
this by first describing an algorithm which may overshoot $\alpha$,
and then making various improvements to ensure that we halt right at
$\alpha$.
 
Since $\alpha$ is recursive, there is a Turing machine which can
enumerate a real $x$ coding a relation on $\omega$ with order type
$\alpha$. Therefore, there is a machine which can write $x$ on the tape
in $\omega$ many steps.  By using every other cell of the scratch tape,
or some such contrivance, we may assume that there is still plenty of
room for further computation without overwriting this real $x$. Now we
will gradually erase elements from the field of the relation coded by
$x$, using the same idea as in the Count-Through Theorem \ref{twotwo}.
In $\omega$ many steps, we can find the particular natural number $n$
which is the least element in the relation coded by $x$, and then, in
another $\omega$ many steps, we can erase all mention of this element
from the field of the relation. Systematically iterating this, we will
successively erase the least element from the relation until there is
nothing left in the relation. Then, with one further $\omega$ limit, we
will know, while futilely searching for the `next' least element, that the
relation is now empty, and we will halt. This algorithm will take
$\omega+(\omega+\omega)\alpha+\omega$ many steps, since we had the
first $\omega$ many steps to write the original real, and then for each
element of the field of that relation, of which there are $\alpha$
many, we had $\omega$ many steps to find that it was the least element,
and another $\omega$ many steps to erase it from the field of the
relation. And then we had $\omega$ many steps on top to discover that
we were actually done (one needs to put the most important flag,
signaling that we are done, on the very first cell of the tape, in
order that the machine will not need extra steps after that last
limit).
 
Let us now improve the algorithm and reduce the number of steps
required.  First, we can combine the initial writing of the relation
with the searching for the least element of the relation, and perform
these functions simultaneously in the first $\omega$ limit. This
produces an algorithm which halts after $(\omega+\omega)\alpha+\omega$
many steps. Second, we can combine the operation of searching for the
next least element and erasing mention of the previous least element
into one $\omega$ limit, by performing them simultaneously. This
produces an algorithm which halts after $\omega\alpha+\omega$ many
steps. Now, we will do infinitely many things simultaneously. Rather
than searching merely for the least element of the relation, we will
search for the least $\omega$ many elements of the relation. This can
be done by systematically searching through the relation, and writing
all the natural numbers from the field which have been encountered so
far on the scratch tape in their order with respect to the relation.
Since the relation is a well order, the initial segments of these
guesses will eventually stabilize to the least $\omega$ many elements
of the relation, and so in the limit we will have written the least
$\omega$ many elements of the order, the others having been `pushed off
to infinity', so to say.  Thus, in one $\omega$ limit we can
simultaneously guess the next $\omega$ least elements of the order,
even while we are erasing the previous least $\omega$ many elements
from the order. This produces an algorithm which will halt after
$\alpha+\omega$ many steps. Remember, the $\omega$ on top is the result
of not being able to recognize in a limit that we have actually erased
everything from the field of the relation. But the following
contrivance will allow us to recognize in a limit that the field of the
relation has already been erased. At any point in our algorithm, the
real coding the relation has some smallest pair $\<n,k>$ (in the natural
ordering of $\N$) which has not been erased. Every time we erase what
is at that time the smallest pair, let the machine flash a flag on and
then off again on the very first cell. In a limit, if there is still
something in the field, then there will be a smallest thing in the
field which survived through the limit, and so this flag will be $0$.
If in a limit the field of the relation is empty, then infinitely often
the smallest element was erased, and so this flag will be $1$. Thus, we
need not go searching for any elements of the field, since by checking
this flag the algorithm can know if there are any elements left in the
field. With this additional procedure, the algorithm will halt
in exactly $\alpha$ many steps.  The ordinal $\alpha$, therefore, is
clockable. \eth
 
Thus the clockable ordinals extend at least up to $\wCK $, the 
supremum of the recursive ordinals. The next theorem shows, however, that
they extend even further than this. 
 
\newtheorem{moreclock}[recclock]{Theorem}
\begin{moreclock}
The ordinal $\wCK +\omega$ is clockable.
\end{moreclock}
 
\proof By \cite{Fef} there is an r.e. relation coding a linear order whose
well-founded part is exactly $\wCK$. Consider the algorithm
which in the first $\omega$ many steps writes this relation on the scratch
tape, and then begins counting through it using the algorithm of 
the previous theorem. At stage $\wCK $, the well-founded part of the
relation has been completely erased. Thus, in another $\omega$ many steps, 
the machine will discover that there is no next least element. The algorithm
therefore can halt at stage $\wCK +\omega$.\eth
 
We invite the reader to ponder the curious question whether the ordinal
$\wCK $ itself is clockable; this we will answer later. For
now, however, we hope to tease the reader by refusing to state whether
the Gap Existence Theorem, following the Speed-up Lemma, hints at the
answer to this question or distracts one from it.
 
\newtheorem{clockiff}[recclock]{Speed-up Lemma}
\begin{clockiff}
If $\alpha+n$ is clockable for some natural number $n$, then $\alpha$ is 
clockable. 
\label{speedup}
\end{clockiff}
 
\proof  Certainly the converse holds; by adding extra states one
can always make a computation take a certain finite number of steps
longer.  The hard part is to make a computation shorter. So, suppose
that $\alpha+n$ is clockable by the program $p$. By the previous
remarks we may assume that $\alpha$ is a limit ordinal. Consider the
operation of $p$ on input 0.  At stage $\alpha$, the first $n$ cells of
the tapes are set up in such a way, say $\vec a$, that the machine
will halt in exactly $n$ additional steps. What we will do is design a program
which will be able to foresee that the first $n$ cells are like that.
And this is how. We will simulate the operation of $p$ on input 0 on
every other cell of the scratch tape in such a way that $\omega$ many
steps of actual computation will simulate $\omega$ many steps of
computation. In the remaining space, we will flash a master flag, on
the first cell, every time one of the first $n$ simulated cells is not $0$
when the corresponding cell in $\vec a$ is $0$.  Thus, in a limit, the
master flag is $0$ only when the first simulated cells have $0$s where
$\vec a$ has $0$s. Secondly, every time each of the first $n$
simulated cells has gotten a $1$, when the corresponding cell in $\vec a$ has
a $1$, we flash a secondary flag (i.e. we flash the
flag after all of the cells have done it since the last flash). Thus,
in a limit, the secondary flag is $1$ when all of the first cells have
a $1$ which should have a $1$ according to $\vec a$. Thus, in a limit,
the first $n$ cells are set up to make the simulation halt when the master
flag is off, and the secondary flag is on. And this can be checked
right at the limit, since the head is actually looking at three cells at
once---one on each of the tapes. So $\alpha$ is clockable.\eth
 
Any child who can count to $89$ can also count to $63$; counting to a
smaller number is generally considered to be easier. Could this fail
for infinite time Turing machines?  Could there be gaps in the
clockable ordinals---ordinals to which infinite time machines cannot
count, though they can count higher? The answer, surprisingly, is Yes.
 
\newtheorem{existgaps}[recclock]{Gap Existence Theorem}
\begin{existgaps}
There are gaps in the clockable ordinals. In fact, the first gap above 
any clockable ordinal has size $\omega$. 
\end{existgaps}
 
\proof  Suppose $\alpha$ is clockable and $\beta$ is the least
non-clockable ordinal above $\alpha$. The ordinal $\beta$ must be a
limit ordinal, and, by the Speed-up Lemma \ref{speedup}, there are no clockable
ordinals between $\beta$ and $\beta+\omega$.  Let us now show that
$\beta+\omega$ is clockable. This will imply that the first gap beyond $\alpha$
has size $\omega$. In order to halt at $\beta+\omega$, our basic strategy 
will be to recognize $\beta$ as the least stage beyond $\alpha$ at which
no infinite time algorithm halts. This recognition will take an extra $\omega$
many steps, and so $\beta+\omega$ will be clockable.  Consider the
algorithm which simulates the computations of $\varphi_p(0)$ for every 
program $p$. That is, by the contrivance of thinking of the
scratch tape and output tape as divided into $\omega$ many scratch
tapes and output tapes, we will simulate, for every program $p$, the
computation of $p$ on input $0$. By setting things up properly, we can
arrange that for every $\omega$ steps of our computation, $\omega$ many
steps are performed in each of the simulated computations.  Since
$\alpha$ is clockable, one of these simulations, for some fixed program
$p_0$, takes $\alpha$ many steps. After waiting for this program to
halt in the simulation, our algorithm will keep careful track of when
the simulated programs halt. When a stage is found for which none of
the simulations halt, then we have found $\beta$ and we halt. Of
course, it took us an extra $\omega$ many steps to recognize that none
of the simulated computations halted at that stage, so our computation
takes $\beta+\omega$ many steps. Thus, $\beta+\omega$ is clockable.\eth
 
These gaps are a bit mysterious. The following lemmas
reveal a little of their structure.
 
\newtheorem{biggaps}[recclock]{Big Gaps Theorem}
\begin{biggaps}
The gaps in the clockable ordinals become large. Indeed, for every
clockable ordinal $\alpha$, there are gaps of size at least $\alpha$ in
the clockable ordinals.
\end{biggaps}
 
\proof  Assume without loss of generality that $\alpha$ is a 
limit ordinal.
What we need to prove is that there are limit ordinals
$\beta$ such that no ordinal between $\beta$ and $\beta+\alpha$ is
clockable, though there are ordinals beyond this that are clockable.
And this is how we will do it. Consider the algorithm which searches
for a gap of size $\alpha$. It does this by simulating every program
$p$ on input $0$ and keeping track of which programs have halted.
Whenever it finds a stage at which none of the programs have halted, it
starts clocking $\alpha$ on the side. For each step of the simulation,
it runs one step of the $\alpha$-clock, and pays attention to determine
whether the clock or the gap runs out first. If the gap runs out first,
then the clock is reset, and the machine searches for the next gap. If the
$\alpha$-clock runs out first, revealing a gap of size $\alpha$, then
the machine halts. If there are no gaps of size $\alpha$, then this
machine will run through all the clockable ordinals, until it finds the
first $\alpha$ many non-clockable ordinals above all of the clockable
ordinals, and then halts.  This is a contradiction since it produces a
clockable ordinal above all the clockable ordinals. Hence, gaps of size
$\alpha$ must exist lower down.\eth
 
\newtheorem{wrogaps}[recclock]{Many Gaps Theorem}
\begin{wrogaps}
There are many gaps in the clockable ordinals. Indeed, if $\alpha$ is a
writable ordinal, then there are at least $\alpha$ many gaps of size at
least $\alpha$ in the clockable ordinals. Moreover, if $\alpha$ is
either clockable or writable, then the exact ordinal number of gaps of size at
least $\alpha$ is neither clockable nor writable.
\end{wrogaps}
 
\proof  We may again assume that $\alpha$ is a limit ordinal.
Recall that a writable ordinal is one which is the order type
of a relation coded by a writable real. The key idea of this argument
is that clockable and writable ordinals both provide a sort of clock by
which to measure the length of time elapsed during a simulated
computation. That is, we can keep track of how long a simulation takes
either by counting through a relation coding $\alpha$, in the case that
$\alpha$ is writable, or by performing one additional step of
computation in an $\alpha$-clock, in the case that $\alpha$ is
clockable. So, suppose now that $\alpha$ is either clockable or
writable, and that the exact number of gaps of size at least $\alpha$ is
$\beta$, where $\beta$ is either clockable or writable. Consider the
algorithm which simulates the computation of all programs $p$ on input
$0$, searching for gaps. Each time it finds a gap, it counts, using the
$\alpha$-clock, to see if the gap has size at least $\alpha$. If so,
then the algorithm counts once on the $\beta$-clock. When the
$\beta$-clock runs out, then all the gaps have been counted through
(and consequently the algorithm has computed beyond every clockable
ordinal). By halting when this occurs, the algorithm halts beyond all
the clockable ordinals, a contradiction.  Thus, $\beta$ must be
neither clockable nor writable.\eth
 
\newtheorem{initseg}[recclock]{No Gaps Theorem}
\begin{initseg}
There are no gaps in the writable ordinals. 
\label{fourone}
\end{initseg}
 
\proof  Suppose $\alpha$ is a writable ordinal. Thus, there is a
program $p$ which on input $0$ writes a real $x$ coding a relation
$\prec_x$ with order-type $\alpha$. If $\beta<\alpha$, then there is
some natural number $n$ which is the $\beta^{\rm th}$ element of the
order $\prec_x$. We can now direct a machine to first write $x$ on the
tape, and then delete from the field of the relation every element
which is not below $n$ with respect to $\prec_x$. After this, the
machine has written ${\prec_x}\restrict n$, which has order type $\beta$,
and so $\beta$ is writable.\eth
 
Before proving the next theorem, let us define a real $x$ to be
{\it accidentally} writable when it appears on one of the tapes during a
computation, but not necessarily as the output of a computation.
Similarly, a real $x$ is {\it eventually} writable 
when there is a nonhalting 
infinite time computation, on input 0, which eventually writes x on 
the output tape; that is, beyond some stage, the real appearing on 
the output tape is x.
Thus, it is clear that every writable real is
eventually writable, and every eventually writable real is accidentally
writable.
 
\newtheorem{otco}[recclock]{Order-type Theorem}
\begin{otco}
The classes of clockable and writable ordinals have the same order-type. 
And this order-type is equal to the supremum of the writable ordinals, which
is neither clockable nor writable, though it is eventually writable.
\label{threetwelve}
\end{otco}
 
\proof First we will show that if $\beta$ is clockable, then the order-type of
the clockable ordinals up to $\beta$ is writable. We will design a supertask
algorithm which will write the relation in which $p$ is less than $q$
when $p$ halts before $q$ on input $0$ and both halt in less than
$\beta$ many steps. Actually, this relation is a pre-wellorder---two
programs are equivalent when they halt at the same stage---and so we actually
want to include into the field of the relation only the least element
of each equivalence class. This produces a relation $\triangleleft$
whose order-type is the same as the order-type of the clockable
ordinals below $\beta$. To write the relation $\triangleleft$, we
simply simulate the operation of all programs on input $0$,
while simultaneously running a clock for $\beta$. In each block of 
$\omega$ many steps of computation we simulate one step 
of computation for each of the programs.
At the beginning and also after each limit stage, we also 
compute one step of the $\beta$-clock.  By keeping careful track
of which programs have halted, we can gradually write the relation
$\triangleleft$ on the output tape. Thus, if program $p$ halts, and
then $q$ halts, but the $\beta$-clock is still running, then we will
ensure that $p\triangleleft q$ on the output tape (while making sure to
include only the least element of each equivalence class in the
field).  When the $\beta$ clock halts, then we halt, and we have
written $\triangleleft$, as desired. It follows from this that the
order-type of the writable ordinals is at least that of the clockable
ordinals.
 
Now we must show the converse.
For the remainder of the paper, let 
$\lambda$ be the supremum of the
writable ordinals. Since there are no gaps in the writable ordinals,
$\lambda$ is also the least non-writable ordinal. We must show that
there are at least $\lambda$ many clockable ordinals. For this, it
suffices to show that whenever $\alpha$ is writable, then there are at
least $\alpha$ many clockable ordinals. Consider the algorithm which
first writes a real coding $\alpha$ on the tape, and then begins
counting through $\alpha$, gradually erasing the field of the relation
in order.  For each natural number $n$, we could design a program which
halts when $n$ is the least element of the field of the relation. That
is, for each $n$, we could arrange for a machine $p_n$ to halt when it has
erased the relation up to $n$. With different values of $n$, these
programs will halt at different times, and the order-type of their various
halting times will be $\alpha$. So there are at least $\alpha$ many
clockable ordinals, and the first part of the theorem is proved.
 
For the second part, we have already mentioned that $\lambda$ is not
writable. Suppose it were clockable. Then we could simulate all the
programs $p$ on input $0$, keeping track of which computations have
halted.  Every time we find a stage at which one of the computations
halts, we run one step of the computation which clocks $\lambda$. When
the clock runs out, we halt. Since there are exactly $\lambda$ many
clockable ordinals, this algorithm will halt beyond all the
clockable ordinals, a contradiction.
 
It remains only to show that $\lambda$ is eventually writable. But this
is easy. As in the argument above, define that $p\triangleleft q$ when
$\varphi_p(0)$ halts before $\varphi_q(0)$,
where again we include only the least program from each equivalence class of 
all programs halting at the same time.
This relation is eventually
writable by the algorithm which simply simulates all computations
$\varphi_p(0)$, and outputs Yes whenever it determines that one
program halts before another. Since this
relation has order-type $\lambda$, we have proved that $\lambda$ is
eventually writable.\eth
 
Let us now prove a few closure theorems for the clockable ordinals.
 
\newtheorem{clockadd}[recclock]{Theorem}
\begin{clockadd}
If $\alpha$ and $\beta$ are clockable, so are $\alpha+\beta$
and $\alpha\cdot\beta$.
\end{clockadd}
 
\proof  For addition, it suffices to consider only the case when
$\beta$ is at least $\omega^2$. This will allow us to be a bit sloppy.
Consider the supertask in which we first clock $\alpha$, and then erase the
tape completely, and then clock $\beta$. This takes
$\alpha+\omega+\beta$ many steps, but by our assumption on $\beta$ it
follows that $\omega+\beta=\beta$, and consequently our algorithm took
$\alpha+\beta$ many steps.
 
The basic idea for ordinal multiplication is that we can clock through
$\beta$, and for each tick of that clock, we clock $\alpha$. This will
take $\alpha\cdot\beta$ many steps.\eth
 
We can also handle large sums of clockable ordinals, of the form 
$\Sigma \alpha_p$, where the sum is taken over a computable set of 
programs $p$, each of which clocks the corresponding ordinal $\alpha_p$. 
That is, if there is some writable real $x$ such that the sum only includes 
the programs $p$ such that $x(p)=1$, we refer to the sum 
$\Sigma\alpha_p$ as a {\it computable sum}. 
 
\newtheorem{closure}[recclock]{Theorem}
\begin{closure}
The supremum of the clockable ordinals is closed under 
infinite time computable addition. In this sense, it is
supertask inaccessible.
\end{closure}
 
\proof For the remainder of 
the paper, let $\gamma$ be the
supremum of the clockable ordinals.
Suppose $x$ is a writable real coding a sequence of programs $p$ which
clock some ordinals $\alpha_p$. The intended sum is
$\Sigma_{x(p)=1}\alpha_p$.  Consider the supertask algorithm which first writes
$x$ and then runs in turn each program $p$ such that $x(p)=1$, and then
halts. This algorithm takes at least $\Sigma_{x(p)=1}\alpha_p$ many 
steps, so $\gamma$ is closed under computable addition.\eth
 
The next theorem shows that there are long stretches of clockable ordinals
without any gaps. 
 
\newtheorem{gaplessblocks}[recclock]{Gapless Blocks Theorem}
\begin{gaplessblocks} 
There are large gapless blocks of clockable ordinals. Indeed, 
if $\alpha$ is writable in $\delta$ many steps, then $\delta+\beta$ is clockable
for any $\beta\leq\alpha$.
\end{gaplessblocks}
 
\proof  Suppose $\alpha$ is writable, so that there is a program
which writes a real $x$ coding a relation with order-type $\alpha$, in
$\delta$ many steps. Now suppose $\beta \leq \alpha$. It suffices to treat
the case when $\beta$ is a limit ordinal. There
must be some natural number $n$ which is the $\beta^{\rm th}$ element
in the relation coded by $x$.  Consider the algorithm which first
writes $x$, and then counts through the relation coded by $x$, ignoring (and
erasing) any part of the field of the relation at $n$ or above. This will 
take an additional $\beta$ many steps. By the technique of the Recursive 
Clocks Theorem \ref{threeone}, the algorithm can recognize in a limit that the relation 
has already been completely erased. Thus, $\delta+\beta$ is clockable.\eth
 
\newtheorem{woclocka}[recclock]{Theorem}
\begin{woclocka}
If $\alpha$ is either clockable or writable, then the set of reals coding 
well orders of length less than $\alpha$ is decidable.
\end{woclocka}
 
\proof  Given an input $x$, we count through it, using the algorithm of
Theorem~\ref{twotwo}, while at the same time counting through $\alpha$,
either with a clock, or else by first writing a real coding $\alpha$
and gradually erasing elements from the field, depending on whether
$\alpha$ is clockable or writable. By paying attention to which
computation runs out first, we can know if the given real codes a well
order with length less than $\alpha$.\eth
 
The following basic question remains open. At the end of the paper we will 
identify necessary and sufficient conditions for an affirmative answer.
 
\newtheorem{Que}[recclock]{Question}
\begin{Que}
Is every clockable ordinal writable?
\end{Que}
 
\section{The Infinite Time Halting Problems}
 
The halting problem, ubiquitous in classical computability theory, has
an infinite time analog which we will analyze in this section.
Officially, we define the halting problems as the sets
$H=\set{(p,x)\mid p \hbox{\rm\ halts on input }x}$, and, the light-face
version, $h=\set{p\mid p \hbox{\rm\ halts on input }0}$. We will see
later that, unlike their classical analogs, these two sets are not
equivalent.
 
\newtheorem{exacthalt}{Halting Problem Theorem}[section]
\begin{exacthalt}
The halting problems $h$ and $H$ are semi-decidable but not decidable.
\label{haltprob}
\end{exacthalt}
 
\proof  Clearly the halting problems are semi-decidable: to determine
if $\varphi_p(x)$ halts, one simply simulates the computation, and if
it ever does halt, output Yes; otherwise keep simulating.  So both $h$
and $H$ are semi-decidable.
 
Let us now prove that they are not decidable. Suppose the halting
problem $H$ was decided by the infinite time computable function $r$.
Thus, $r(p,x)=1$ when $\varphi_p(x)\converges$, and otherwise
$r(p,x)=0$.  Let $q$ be the program which computes the following
function:  $$\varphi_q(p)=\cases{\uparrow& r(p,p)=1\cr 1&
r(p,p)=0\cr}$$ Now simply observe that $q$ halts on input $q$ iff $q$
does not halt on input $q$, a contradiction. So $H$ is not decidable.
 
To see that $h$ is not decidable, we will use the Recursion Theorem.
Suppose $h$ was decided by some computable function $r$.  By the
Recursion Theorem \ref{Recthm}, there is a program $q$ such that
$$\varphi_q(n)=\cases{\uparrow& r(q)=1\cr 1& r(q)=0}$$ Thus, taking
$n=0$, we see that the program $q$ halts on input $0$ exactly when the
program $q$ does not halt on input $0$, a contradiction.\eth
 
We will sometimes want to refer to approximations to the halting
problem, and so we define $H_{\alpha}=\set{(p,x)\mid p\hbox{\rm\ halts
on input $x$ in fewer than $\alpha$ steps}}$, and $h_{\alpha} =
\set{p\mid p\hbox{\rm\ halts on input $0$ in fewer than $\alpha$
steps}}$.  It is clear from the definition that if $\alpha < \beta$
then $H_{\alpha} \subseteq H_{\beta}$. What is more, as we will now
prove, this inclusion is strict.
 
\newtheorem{H2}[exacthalt]{Possible Lengths Theorem}
\begin{H2}
Infinite time computations come in all possible lengths. Indeed, there
is a single program which, on various input, performs a halting
computation taking any specified non-zero countable ordinal length of time.
\end{H2}
 
\proof Consider the following supertask algorithm. On input $x$, if the
first digit of $x$ is $1$, then the algorithm searches for the next $1$ and 
halts upon finding it. That will take care of computations of non-zero 
finite length.
Now, if the first digit of $x$ is $0$, and the second digit is $1$, then 
the remaining input is interpreted as a limit ordinal to be counted through
according to the algorithm of Theorem \ref{threeone}. That will take 
care of the computations of limit ordinal length. Finally, if the 
first two digits of $x$ are both $0$, then the next $1$ is searched for and 
then moved two spaces to the left, 
while what comes after is interpreted as coding a relation to 
be counted through according to the algorithm of Theorem \ref{threeone}
(leaving that first $1$ in place). 
After this, the algorithm counts through the empty space until it finds again
that first $1$ and then halts. For any limit ordinal $\beta$ and any natural 
number $n$, this kind of input can be arranged to take exactly $\beta+n$ many 
steps. Thus, on various input, this program computes for any desired non-zero
length of time.\eth
 
\newtheorem{H1}[exacthalt]{Corollary}
\begin{H1}
If $\alpha < \beta$ then $H_{\alpha}$ is a proper
subset of $H_{\beta}$.
\end{H1}
 
\newtheorem{H3}[exacthalt]{Theorem}
\begin{H3}
For any limit ordinal $\alpha$, neither $H_\alpha$ nor $h_{\alpha}$ is
$\alpha$-decidable. But if $\alpha$ is clockable, then both $H_{\alpha}$ and
$h_\alpha$ are $\alpha$-semi-decidable and $(\alpha+1)$-decidable.
\end{H3}
 
\proof  For the first part of this theorem, simply observe that if
the function $r$ in the Halting Problem Theorem \ref{haltprob} is
computable in fewer than $\alpha$ many steps, so is the function
computed by the program $q$ which we defined there.  So the
contradiction of that argument goes through.
 
For the second part, assume that $\alpha$ is a clockable limit ordinal,
and consider the algorithm which, on input $(p,x)$ runs $p$ on input
$x$, while at the same time running a program which clocks $\alpha$. If
$p$ halts on $x$ before $\alpha$ is clocked, our machine halts with
output $1$. This occurs before $\alpha$.  Otherwise, our machine never
halts. Thus, $H_{\alpha}$ is $\alpha$-semi-decidable. It follows that
$h_\alpha$ is also $\alpha$-semi-decidable.
 
Lastly, consider the algorithm which operates as in the above
paragraph, except that if the $\alpha$-clock halts, this machine halts
with an output of $0$. The program can be arranged so as to recognize
that the $\alpha$-clock has halted in the $\alpha^{\rm th}$ step, 
thus deciding $H_{\alpha}$, and also $h_\alpha$, in fewer than $\alpha+1$ many
steps.\eth
 
\newtheorem{H5}[exacthalt]{Theorem}
\begin{H5}
If $\alpha$ is writable or clockable, then $H_{\alpha}$ and $h_\alpha$ are 
decidable.
\label{seven}
\end{H5}
 
\proof The previous argument handles the case when $\alpha$ is
clockable.  So assume now that $\alpha$ is writable, and consider the
machine which, upon input $(p,x)$, writes a code for $\alpha$ on the
output tape, and then simulates $p$ on input $x$, erasing an element
from the coding of $\alpha$ after each step in the program.  If $p$
halts on $x$ while there is still a non-empty well-order coded, our
machine halts with a $1$.  Otherwise, if the well-order is completely
erased, our machine halts with a $0$. This clearly allows the machine to
decide membership in $H_\alpha$. It follows that $h_\alpha$ is also
decidable.\eth
 
\newtheorem{H6}[exacthalt]{Theorem}
\begin{H6}
The set $h_\alpha$ is decidable for every 
$\alpha$ below the supremum of the clockable ordinals.
\end{H6}
 
\proof Suppose $\alpha<\beta$, where $\beta$ is the least clockable ordinal
above $\alpha$. Since no computations on input $0$ can halt between $\alpha$
and $\beta$, it follows that $h_\alpha=h_\beta$. And we know $h_\beta$ is
decidable by the previous theorem.\eth
 
\newtheorem{H7}[exacthalt]{Theorem}
\begin{H7}
Let $\gamma$ is the supremum of the clockable ordinals; then 
$H_{\gamma}$ is semi-decidable but not decidable.
\end{H7}
 
\proof  Clearly $H_\gamma$ is semi-decidable, since on input 
$(p,x)$, one simulates the program $p$ on input $x$, while searching 
for a clockable ordinal which is larger. That is, while simulating 
$\varphi_p(x)$, the program also simulates $\varphi_q(0)$ for all programs
$q$, and pays attention to when the simulations halt. If one of the computations
$\varphi_q(0)$ halts after the computation $\varphi_p(x)$, then, and only 
then, may it be concluded that $(p,x)$ is in $H_\gamma$. So $H_\gamma$ is 
semi-decidable. It cannot be decidable, because $h=h_\gamma$ appears as its
$0^{\rm th}$ slice.\eth
 
We foreshadowed the next few theorems at the end of the second section.
If $A$ is a subset of the plane $\R\cross\R$, then the slices of $A$
are the sets $A_y=\set{z\st (y,z)\in A}$. Recall that a settled
snapshot sequence for $\varphi_p(x)$ is a real that codes a
well-ordered sequence of snapshots of the computation of $\varphi_p(x)$
whose last snapshot is either a halting snapshot or else is the
first repeating snapshot (in the strong sense explained
just after Theorem \ref{onetwo}).
 
\newtheorem{uniform}[exacthalt]{No Uniformization Theorem}
\begin{uniform}
There is a decidable subset of the plane $\R\cross\R$, 
all of whose slices are 
nonempty, which does not contain the graph of any computable total function.
\end{uniform}
 
\proof The following set will do:
$$A=\set{(\<p,x>,z)\st \hbox{$z$ codes a settled snapshot sequence 
for $\varphi_p(x)$}}.$$
With a suitable pairing function, we may assume every real has the form
$\<p,x>$, and, since every computation $\varphi_p(x)$ has a settled
snapshot sequence, every slice of $A$ is non-empty.  The set $A$ is
certainly decidable, since to verify that a proposed real does really
code a settled computation sequence is very easy; the machine must
merely check that the successor steps of the computation are modeled
correctly, and that the limit steps are computed correctly from the
earlier snapshots, and that the last snapshot is either in a halting
state or is the first time that the snapshots repeat.  But there can be
no computable function $f$ such that $f(\<p,x>)$ always gives a settled
snapshot sequence for $\varphi_p(x)$, since such a function could be
easily used to decide $H$: one would simply compute $f(\<p,x>)$
and observe if the last snapshot was in a halt state.\eth
 
We should mention that if a decidable subset $A$ of the plane has 
an accidentally writable real in each section, or indeed, merely has for 
every real $y$ an accidentally $y$-writable real in the section $A_y$, 
then $A$ will contain the graph of a computable function. On input $y$, 
simply search for an accidentally $y$-writable real $z$ such that $(y,z)\in A$.
Such a real will eventually be found, and so this algorithm gives a 
computable total function uniformizing $A$. This shows that 
some computations $\varphi_p(x)$ cannot have even accidentally writable 
settled snapshot sequences. 
 
Like the previous theorem, the next identifies a surprising divergence 
from the classical theory. The real $c$ in the theorem is like a 
forgotten melody that you cannot produce on your 
own but which you can recognize when someone sings it to you. 
 
\newtheorem{singleton}[exacthalt]{Lost Melody Theorem}
\begin{singleton}
There is a real, $c$, such that $\set{c}$ is decidable, but $c$ is not writable.
Consequently, there is a constant, total function which is not
computable, but whose graph is nevertheless decidable: $f(x)=c$. 
\label{melody}
\end{singleton}
 
\proof The repeat-point of a computation is the ordinal stage by which
it either halts or repeats. Let $\delta$ be the supremum of the
repeat-points of the computations of the form $\varphi_p(0)$. Thus,
$\delta$ is a countable ordinal in $L$. 
Consequently, by a simple bootstrap argument, there is some smallest 
$\beta\geq\delta$ such that $L_{\beta+1}\models\beta$ is 
countable (one can take $\beta=\sup\beta_n$ where $\beta_0=\delta$ 
and $\beta_n$ is countable first in $L_{\beta_{n+1}}$). 
Thus, since $L_{\beta+1}$
has a canonical well-ordering, there is some real $c\in L_{\beta+1}$
which is least with respect to the canonical $L$ order, such that $c$
codes $\beta$. This is our real $c$. First, we will argue that $c$ is
not writable. Indeed, it is not even accidentally writable. If it were,
then we could solve the halting problem $h$ by searching for an
accidentally writable real that codes an ordinal large enough to see
the repeat-point of the computation in question. Since $c$ codes
$\beta$, which is as large as $\delta$, the real $c$ is large enough,
and so our algorithm will succeed.  This contradicts the fact that $h$
is not decidable. Second, we will argue that $\set{c}$ is decidable.
Given a real $z$, we must decide if $z=c$ or not.  First, we can
determine whether $z$ codes an ordinal or not. Suppose that $z$ codes
$\alpha$.  Next, we can determine whether this ordinal is larger than
$\delta$ by simulating every computation $\varphi_p(0)$ along the order
given by $z$, and determining if we reach the repeat-point before
running out of room. Now comes the complicated part. A hereditarily
countable set $a$ may be coded with a real by first computing the
transitive closure of $\{a\}$, and then, since this is countable, by
finding a relation $E$ on $\omega$ such that
$\<\TC(\{a\}),{\in}>\iso\<\omega,E>$. We can then code the relation $E$
with a real in the usual manner (this coding technique is used
extensively in section $8$). This way of coding sets with reals works
well with infinite time Turing machines. In particular, whether a real
is a code is a $\Pi^1_1$ property; whether two codes code the same set
is a $\Sigma^1_1$ property, as is whether one code codes a set which is
an element of the set coded by another code. These elementary
properties are therefore infinite time decidable.  Indeed, the truth of
any $\Delta_0$ set-theoretic property is decidable in the codes. Now,
using the real $z$ which codes $\alpha$ to organize our construction,
we can write a code for the set $L_{\alpha+1}$ by mimicking the
construction of the $L$ hierarchy along the well order given by $z$.
That is, we first use $z$ to reserve infinitely much room for each
$\beta<\alpha$, and given the code for $L_\beta$ we then write down,
from the definition, the code for $L_{\beta+1}$ in the space we had
reserved.  At limit stages, we simply write down the code for the
union of sets whose codes we have already written down. Using the code
for $L_{\alpha+1}$ that we have thus produced, we can check whether $z$
really is the least code in $L_{\alpha+1}$ for $\alpha$. And we can
check whether $\alpha$ really is the least ordinal above $\delta$ such
that $\alpha$ is countable in $L_{\alpha+1}$ (i.e. that $\alpha$ is
$\beta$).  If $z$ passes all of these tests, then $z$ must be $c$,
otherwise it is not.  So $\set{c}$ is decidable. It follows, for the
second part of the theorem, that the set $\set{(x,y)\st y=c}$ is
decidable; this is the graph of the constant, total function
$f(x)=c$. Thus, the graph of $f$ is decidable, but, since $c$ is not
writable, $f$ is not computable.\eth
 
We know how to code ordinals with reals. The previous proof, however,
shows how to associate a {\it unique} such code to the ordinals which
are countable in $L$. Namely, given any ordinal $\alpha$ which is
countable in $L$, let $\beta$ be the least ordinal above $\alpha$ such
that $L_{\beta+1}$ knows that $\beta$ is countable, and let $x$ be the
least real, in the $L$ order, which codes $\beta$. The ordinal $\alpha$
is the $n^{\rm th}$ element in the order given by $x$, for some $n$,
and so $\alpha$ may be coded with the pair $\<x,n>$. This code is
unique since $\alpha$ determines $\beta$, which determines $x$, and
then $x$ and $\alpha$ determine $n$.  Furthermore, the set of such
codes, as we proved in the theorem, is a decidable set.
 
\section{Oracles}
 
Since infinite time Turing machines naturally compute functions on the
reals, we are pushed towards two distinct kinds of oracles: individual
reals and sets of reals. An individual real $x$ can be used as an
oracle much as it is in the Turing machine context, by adding a special
oracle tape on which $x$ is written out.  This amounts, in effect, to
having another input tape and using $x$ as an additional argument.
Nevertheless, to highlight the oracle nature of the computations, if
$p$ is a program using an oracle tape, then we will denote by
$\varphi_p^x$ the resulting function which uses the real $x$ as an
oracle.  We refer to such functions as the infinite time $x$-computable
functions.
 
Oracles, though, are properly the same type of object as decidable or
semi-decidable sets; in the context of infinite time computations, this
means that we want somehow to use a {\it set} of reals $A$ as an
oracle. Clearly we cannot expect always to be able to write such an
object out on a tape; but somehow we want the machine to be able to ask
membership queries of $A$.  Thus, we propose to add a special oracle
tape, initially filled with zeros at the beginning of a computation,
and to allow the machine during a computation to write on this tape and
then, during a computation, by switching to a special {\it oracle
query} state, to receive the answer Yes or No, on the cell beneath the
head, accordingly as the real on the oracle tape is in $A$ or not.
Thus, a machine with oracle $A$ is allowed to know if $y$ is a member of
$A$, for any $y$ which it is able to write on the oracle tape. The
machine can make as many such queries as it likes. We believe that this
notion of oracle for sets of reals is natural and robust, resembling as
it does the notion of constructibility from a predicate, as in the
definition of $L[A]$. We can denote as usual by $\varphi_p^A(x)$ the
resulting function computed by program $p$ with oracle $A$ on input
$x$, and we refer to such functions as the infinite time $A$-computable
functions, or, normally, just as the $A$-computable functions. If $x$
is a real, i.e.  if $x\in 2^\omega$, let $A_x=\set{ (x\restrict
n)\concat\<0,0,\ldots>\st n\in\omega}$ be the corresponding set oracle
of finite approximations to $x$, concatenated with zeros. It is easy to
see that a function is $x$-computable iff it is
$A_x$-computable, because with the real oracle $x$ we can decide
membership in $A_x$, and with the set oracle $A_x$, we can write $x$ on
the tape. In this way real oracles can be thought of as a special case
of set oracles. (Notice that, in view of the Lost Melody Theorem
\ref{melody}, being $x$-computable and being $\{x\}$-computable are
very different things.  In the former, we can see $x$ written out on
the tape; in the latter, we are allowed only to know Yes or No whether
any real we can produce is equal to $x$. But if $x$ is hard to produce,
this is not helpful, since the answer will most likely be No. Indeed,
this is the key significance of the Lost Melody Theorem \ref{melody}.)
 
If $A$ and $B$ are two oracles, we will say that $A$ is infinite time
computable from $B$, written $A\ITleq B$, when the characteristic
function of $A$ is infinite time $B$-computable. This definition makes
sense for both set and real oracles, if we think of reals as subsets of
$\omega$.  Since it is easy to verify that $\ITleq$ is transitive and
reflexive, we also obtain the notion of infinite time degrees:
$A\ITequiv B$ iff $A\ITleq B$ and $B\ITleq A$. This is an equivalence
relation, and we denote the equivalence classes by $[A]_\infty$. Thus, for a
real $x$, we have that $x\ITequiv A_x$.  We will write $A\ITlt B$ when
$A\ITleq B$ and $A\not\ITequiv B$. And of course we extend the notions
of infinite time semi-decidability, clockability, and writability to
the context of oracles in the obvious way, so that, for example, an
ordinal $\alpha$ is $A$-clockable when there is program using
oracle $A$ which halts on input $0$ in exactly $\alpha$ many steps. 
Thus, a real $x$ is $A$-writable if and only if $A_x \ITleq A$ if and 
only if $x \ITleq  A$. One last bit of
notation: given two oracles $A$ and $B$, we write
$A\oplus B$ to mean an oracle which codes, in some canonical manner,
the information contained in $A$ and $B$. Thus, for example, $A\oplus
B$ could be the reals resulting from adding a $0$ digit to the front
of every real in $A$, and a $1$ to those in $B$. Clearly $A\oplus B$ is
the least upper bound of $A$ and $B$ with respect to $\ITleq$.
 
We now define two jump operators, corresponding to the two halting
problems. Suppose that $A$ is an oracle
(either a set or a real oracle). The strong jump of $A$, denoted
by $A^\Jump$, is simply the halting problem relativized to $A$.  That
is, $$A^\Jump= H^A = \set{(p,x)\mid \varphi_p^A(x)\converges}.$$
Secondly, the weak jump of $A$ is the set $$A^\jump= A \oplus
h^A= A \oplus \set{p\mid\varphi_p^A(0)\converges}.$$ It may seem odd that we
explicitly include the factor $A$ into $A^\jump$, but the fact is that
some sets of reals $A$ are sufficiently complex that they are not
computable from any real; in particular, such sets $A$ are not
computable from $h^A$ alone. Since we definitely want $A^\jump$ to
compute $A$, we are led to the definition above. When $A$ is a real,
then it is not difficult to see that $A^\jump\ITequiv h^A$, and
consequently in this case we don't need to explicitly include $A$.
Up to equivalence, $0^\jump$ is just $h$, and $0^\Jump$ is $H$. 
It is easy to see that any $A$-semi-decidable set is computable from 
$A^\Jump$, and any $A$-semi-decidable real is computable from 
$A^\jump$. 
 
\newtheorem{strict}{Jump Theorem}[section]
\begin{strict}
$A \ITlt A^\jump \ITlt A^\Jump$.
\end{strict}
 
\proof  It is easy to see that $A \ITleq A^\jump \ITleq A^\Jump$, since
$A^\jump$ explicitly computes $A$, and $h^A$ appears as the 
$0^{\rm th}$ slice of $A^\Jump$.
What remains is to prove the strict relations. 
The first argument is merely a relativization of Theorem \ref{haltprob}. 
Relativizing that proof, we see that $h^A$ is not computable from $A$, and 
therefore, $A^\jump$ is not computable from $A$, so $A \ITlt A^\jump$. 
 
Let us now prove that $A^\jump\ITlt A^\Jump$. If not, then there is some
program $q$ which computes $A^\Jump$ from $A$ and $z=h^A$. That is, 
$$\varphi_q^{A\oplus z}(p,x)=\cases{1,& if $\varphi_p^A(x)\converges$;\cr
                        0,& if $\varphi_p^A(x){\uparrow}$.\cr}$$

Let $f(r)$ be the program which halts on input $x$, using oracle $A$, 
exactly when $\varphi_q^{A\oplus x}(r,x)=0$. Thus, by the recursion 
theorem, there is a program $r$ such that
$\varphi_r^A(x)\converges$  exactly when $\varphi_q^{A\oplus x}(r,x)=0$. 
But a special case of this is 
$\varphi_r^A(z)\converges\leftrightarrow \varphi_q^{A\oplus z}(r,z)=0$,
which contradicts 
the assumption on $q$.\eth
 
The argument just given actually establishes the following corollary. 
 
\newtheorem{fromareal}[strict]{No Reals Corollary}
\begin{fromareal}
The set $A^\Jump$ is not computable from $A\oplus z$ for any real $z$.
In particular, $0^\Jump$ is not computable from any real.
\end{fromareal}
 
\newtheorem{absorption}[strict]{Absorption Theorem}
\begin{absorption}
$A^{\jump\Jump}\equiv A^\Jump$. Indeed, for any ordinal $\alpha$ which is
$A^\Jump$-writable, $A^{\jump^{(\alpha)}\Jump}=A^\Jump$. 
\end{absorption}
 
\proof Let us prove the first equation first. Since $A\ITleq A^\jump$,
it is clear that $A^\Jump\ITleq A^{\jump\Jump}$. It remains to show the
converse relation. If $p$ is a program, let $f(p)$ be a program such
that $\varphi_{f(p)}^A(x,y)=\varphi_p^{A\oplus y}(x)$. 
We may assume that $f$ is computable. Now simply compute 
$$(p,x)\in A^{\jump\Jump}=(A\oplus
h^A)^\Jump\Iff \varphi_p^{A\oplus h^A}(x)\converges$$
$$\Iff
\varphi_{f(p)}^A(x,h^A)\converges\Iff (f(p),\<x,h^A>)\in
A^{\Jump}.$$ 
Thus, $A^{\jump\Jump}\ITleq A^\Jump$, and we are done.
 
Now, consider the more complicated equation
$A^{\jump^{(\alpha)}\Jump}=A^\Jump$, where $\alpha$ is coded by some
real $z$ which is $A^\Jump$ writable. We define the iterates
$A^{\jump^{(\alpha)}}$ with respect to $z$ by induction on $\alpha$, so
that $A^{\jump^{(\alpha)}}=A\oplus w_\alpha$ for some real $w_\alpha$.
We begin of course with $w_0=0$, so that $A^{\jump^{(0)}}=A\oplus 0$.
At successor stages, we want
$$A^{\jump^{(\beta+1)}}=(A^{\jump^{(\beta)}})^\jump= (A\oplus
w_\beta)\oplus h^{A^{\jump^{(\beta)}}},$$ and so we simply let
$w_{\beta+1}=w_\beta\oplus h^{A^{\jump^{(\beta)}}}$.  At limit stages
$\delta$, we let $w_\delta=\oplus_{\beta<\delta}w_\beta$, using the
real $z$ to organize the information. Now let us prove the theorem by
induction on $\alpha$. Successor stages follow directly from the
previous paragraph, since
$$A^{\jump^{(\beta+1)}\Jump}=A^{\jump^{(\beta)}\jump\Jump}=
A^{\jump^{(\beta)}\Jump}=A^\Jump,$$ where the first equality follows
from the definition of $A^{\jump^{(\beta+1)}}$, the second from the
previous absorption argument, and the third by the induction
hypothesis. Now suppose $\delta$ is a limit ordinal, and the result
holds for every $\beta<\delta$. Since $A^{\jump^{(\delta)}}=A\oplus
w_\delta$, where $w_\delta=\oplus_{\beta<\delta}w_\beta$ (using the
real $z$ to organize the information), it follows that $w_\delta$ is
computable from $A^\Jump$ since the $w_\beta$ are uniformly computable
from $A^\Jump$.  We are really just iterating the argument of the
previous paragraph along the order coded by the real $z$. So, to finish
the theorem, we just have to argue that $(A\oplus w_\delta)^\Jump$ is
computable from $A^\Jump$. To see why this is so, let $q$ be the
program such that $\varphi_q^A(y,p,x)=\varphi_p^{A\oplus y}(x)$. Thus,
$(p,x)\in (A\oplus w_\delta)^\Jump\leftrightarrow (q,\<w_\delta,p,x>)\in
A^\Jump$. And with $A^\Jump$ we can compute this latter property, so
the proof is complete.\eth
 
The next theorem shows that the infinite time jump operators $\jump$
and $\Jump$ jump much higher than the Turing jump, even when the Turing
jump is iterated an enormous number of times.
 
\newtheorem{turing}[strict]{Jump Closure Theorem}
\begin{turing}
Every infinite time degree is closed under the Turing jump operator.
Indeed, for any real $x$ and any writable ordinal $\alpha$, the
$\alpha^{\rm th}$ Turing jump of $x$ is still infinite time equivalent to
$x$.
\end{turing}
 
\proof The Turing jump of a real $x$ is defined to be $x'=\set{e \in
\N: \set{e}^x(e)\converges}$, where $\set{e}^x(e)$ denotes the Turing
machine computation of program $e$ on input $e$ with oracle $x$.  The
$\alpha^{\rm th}$ jump of $x$, relative to a real $z$ coding a
relation on $\omega$ of length $\alpha$, is the subset of the plane
$\omega\cross\omega$ whose $(\beta+1)^{\rm th}$ column is the Turing
jump of the $\beta^{\rm th}$ column, and at limits a sum is taken using
the canonical order given by $z$. The point is that the Turing machine
halting problem is infinite time computable in $\omega$ many steps, and
so an infinite time machine can systematically compute the iterates of
the Turing jumps by simply solving the Turing machine halting problem
for the various oracles, and using the coding of $\alpha$ by the
real $z$ to organize the information. \eth
 
\newtheorem{Heqh}[strict]{Theorem}
\begin{Heqh}
For any oracle $A$, $A^\jump \ITequiv H^A_{\gamma^A}$, where $\gamma^A$ is the
supremum of the $A$-clockable ordinals.
\end{Heqh}
 
\proof  First, let us argue that $H^A_{\gamma^A} \ITleq A^\jump$.
For this, consider the algorithm which on input $(p,x)$ first consults
$A^\jump$ to find which programs $r$ halt with oracle $A$ on input $0$.
Then, the algorithm simulates the computation $\varphi^A_r(0)$ for all
these $r$, at the same time simulating the computation
$\varphi^A_p(x)$. The computations $\varphi_r^A(0)$ act as a clock
which will run out at $\gamma^A$.  If the computation of
$\varphi_p^A(x)$ finishes before all the clocks have run out, then
$(p,x)$ is in $H_{\gamma^A}^A$, and we output Yes. Otherwise, the
clocks all run out first, and we output No.
 
Second, we will argue that $A^\jump\ITleq H^A_{\gamma^A}$. Certainly
$A$ is computable from $H^A_{\gamma^A}$, by considering the algorithm
which halts on elements of $A$ and otherwise does not halt (using $A$
as an oracle). So it remains only to show that $h^A$ is computable from
$H^A_{\gamma^A}$. But this is clear, since a program $p$ with oracle
$A$ on input $0$ will halt before $\gamma^A$ if it halts at all,
because the length of the computation will be an $A$-clockable ordinal.
So $p\in h^A\Iff (p,0)\in H^A_{\gamma^A}$. Thus, both $A$ and $h^A$ are
computable from $H^A_{\gamma^A}$, and so $A^\jump\ITleq
H^A_{\gamma^A}$, as desired.\eth
 
\newtheorem{deltafat}[strict]{Theorem}
\begin{deltafat}
$\Delta^1_2$ is closed under the jump operators $\jump$ and $\Jump$.
\end{deltafat}
 
\proof  We simply have to check that if $A$ is $\Delta^1_2$, then so is
$A^\Jump$. But $(p,x)$ is in $A^\Jump$ exactly when there is a real
coding a well-ordered sequence of snapshots according to program $p$ on
input $x$ with oracle $A$ which shows the computation to halt. Let us
say this in more detail: $(p,x)$ is in $A^\Jump$ exactly when there is
a real $z$ which codes a well-ordered sequence of snapshots such that
first, the initial snapshot is the starting configuration of
$\varphi^A_p (x)$; second, at every step of this sequence either an
ordinary computation was performed by $p$ to obtain the next snapshot,
or, if $z$ indicates that an oracle query was made to which the answer
$z$ provided was Yes, then there is a real $y$ which is the real the
query was made about, and this real is in $A$, or, if $z$ indicates
that an oracle query was made to which the answer $z$ provided was No,
then there is a real $y$ which is the real that the query was made
about, and this real is not in $A$; third, the snapshots at limit
stages are obtained from the previous snapshots according to the
$\limsup$ rule; and fourth, that the final snapshot shows the
computation to have halted. Thus, $A^\Jump$ is $\Sigma^1_2$. Similarly,
$(p,x)$ is in $A^\Jump$ exactly when every real coding a well-ordered
sequence of snapshots according to $p$ on input $x$ with oracle $A$
which is settled shows the computation to halt. This is $\Pi^1_2$,
and so $A^\Jump$ is $\Delta^1_2$, as desired.\eth

Let us now analyze the complexity of the infinite time degree relation. Recall 
that $x\ITleq y$ when $x$ is infinite time computable from $y$. 
 
\newtheorem{degcomp}[strict]{Theorem}
\begin{degcomp}
The relation $x\ITleq y$ is semi-decidable but not decidable.
\label{sixseven}
\end{degcomp}
 
\proof  Notice that $x\ITleq y$ exactly when there is a program
$p$ such that $\varphi_p^y(0)=x$. So we can simply try them all out.
That is, on input $x$ and $y$, simultaneously, for each program $p$,
simulate the computation of $\varphi_p^y(0)$, and, when and if the
simulations halt, check if the output is $x$. If so, output Yes. This
algorithm gives the affirmative answers to $x\ITleq y$, and therefore
that relation is infinite time semi-decidable.
 
Now let us argue that it is not decidable. First note that there are
accidentally writable reals that are not decidable, since $0^\jump$ is
such a real.  Assume towards a contradiction that $\ITleq$ is
decidable, and consider the following supertask algorithm: for each program $p$
simulate the computation of program $p$ on input $0$, and after each step
of the simulated computation, check if the real on the simulated tape
is decidable (this is possible by our assumption). If a nondecidable
real is found, then write it on the output tape and halt. Since there
are accidentally writable reals which are not decidable, this algorithm
will halt with a nondecidable real written on the tape. But this is a
contradiction, since such a real cannot be writable.\eth
 
For the next theorem, denote by $\lambda^A$ the supremum of the
$A$-writable ordinals, which is also the order-type of the
$A$-clockable ordinals. And, as before, let $\gamma^A$ 
be the supremum of the
$A$-clockable ordinals.
 
\newtheorem{jumpup}[strict]{Theorem}
\begin{jumpup}
The ordinal $\lambda^A$ is $A^\jump$-writable in $\gamma^A$ many steps, and 
it is eventually $A$-writable, but not $A$-writable. 
The eventually $A$-writable ordinals are the same as the 
eventually $A^\jump$-writable ordinals.
\end{jumpup}
 
\proof  The first part of this theorem follows by simply relativizing
the argument of \ref{threetwelve} to the context of the oracle $A$.
This algorithm involves computing the relation $p\triangleleft q$ iff
the computation $\varphi_p^A(0)$ halts before $\varphi_q^A(0)$. The
relation has rank $\lambda^A$, and the natural way to write it takes
$\gamma^A$ many steps (one must consult $A^\jump$, asking which
computations $\varphi_p^A(0)$ will halt, in order to know that the
computation has finished).  By omitting this last part, about
consulting $A^\jump$ to know whether the algorithm has finished, 
it follows that $\lambda^A$ is eventually $A$-writable. And 
by relativizing Theorem 3.8 again, we see that $\lambda^A$ is not $A$-writable.

Since $h^A$ is eventually $A$-writable,
any ordinal which is eventually $A^\jump$-writable is
actually eventually $A$-writable, because with only
$A$ as an oracle, we can compute approximations to $h^A$,
and then, with these approximations, run the program
which eventually writes an ordinal from $A^\jump$.
Eventually, this algorithm will have the true $h^A$ written,
and the correct ordinal will be eventually written. \eth
 
\section{The Structure of Infinite Time Degrees}
 
In this section we would like to give an introductory analysis of the
structure of the infinite time degrees. Initially, one might hope to mimic
many of the results from classical recursion theory, perhaps even giving some
priority arguments, but this hope must be tempered by the realization
that one cannot preserve a computation by fixing only finitely much of
an oracle; after all, the computation to be preserved may have used the
entire oracle. Nevertheless, we have managed to prove a few facts, and
content ourselves to leave the vast bulk of results for others to
prove.
 
Let us first prove some elementary facts.
 
\newtheorem{incom}{Theorem}[section]
\begin{incom}
There are incomparable infinite time degrees in the reals.
\end{incom}
 
\proof  This proof also works in the Turing degree case.
Assume that the infinite time degrees are linearly ordered. Since every
initial segment of this order is countable---there being only countably
many programs---it follows that there are at most $\aleph_1$ many
degrees, and consequently the Continuum Hypothesis holds. But the
Continuum Hypothesis fails in a forcing extension. Consequently, in
such extensions there are incomparable degrees. Furthermore, by Theorem
\ref{sixseven}, the assertion that there are incomparable degrees is a
$\Sigma^1_2$ assertion, and thus, by the Shoenfield absoluteness
theorem, it is absolute to any forcing extension. Hence, there must
have been infinite time incomparable reals originally.\eth
 
The previous theorem can be improved to the following.
 
\newtheorem{antichain}[incom]{Antichain Theorem}
\begin{antichain}
There is a countable sequence of reals, no one of which is infinite time 
computable from the others. 
\end{antichain}
 
\proof  Since it is well known that no real in a sequence of mutually
generic Cohen reals can be constructed from the rest of the sequence,
the assertion is true after forcing to add $\omega$ many Cohen reals.
Furthermore, the assertion is $\Sigma^1_2$. Hence, by the Schoenfield
absoluteness theorem, it must have been true originally.\eth
 
The next two theorems are straightforward adaptations of the 
classical arguments (see, e.g., \cite{So}). 
 
\newtheorem{universal}[incom]{Corollary}
\begin{universal}
Any countable partial order embeds into the real infinite time degrees. 
\end{universal}
 
\proof Suppose $\<\omega,\preceq>$ is a partial order on $\omega$, and 
$\<z_n|n\in \omega>$ is a sequence of reals, no one of which is computable 
from the others and $\<\omega,\preceq>$. It is straightforward
to verify that the 
map $n\mapsto\oplus\set{z_i\mid i\preceq n}$ embeds the partial order into 
the real infinite time degrees.\eth
 
\newtheorem{simple}[incom]{Theorem}
\begin{simple}
In the infinite time degrees, a simple set exists. That is, there is a
semi-decidable set $S$ whose complement is infinite but contains no
infinite semi-decidable set.
\end{simple}
 
\proof  For this argument, we can simply mimic Post's proof in the
classical argument. Namely, for each program $p$, start
computing $\varphi_p(n)$ for all the values of $n\in\omega$, and let
$f(p)$ be the first $n>2p$ for which $\varphi_p(n)\converges$. Let
$s=\ran(f)$. Thus, $s$ is a subset of $\omega$, and it is clearly
semi-decidable. Also, the complement $\omega\setminus s$ is infinite
since $s$ contains at most $p$ many elements below $2p$. But $s$ meets
every infinite semi-decidable subset of $\omega$, so the complement
$\omega\setminus s$ contains no infinite semi-decidable set. To
construct a set of reals with the same property, simply observe that we
can regard $\omega\subset \R$, and then let
$S=(\R\setminus\omega)\union s$. Thus, $S\subset \R$ is semi-decidable
with an infinite complement in $\R$ which contains no infinite
semi-decidable set.\eth
 
\newtheorem{incom0}[incom]{Questions{\footnotesize$^2$}}
\begin{incom0}
\bd
\item
\item Is there a semi-decidable set of reals $S$ whose complement has size the continuum 
but contains no infinite semi-decidable set? 
\item Are there incomparable degrees below $0^\jump$? Below $0^\Jump$? 
\item Are there any noncomputable degrees below $0^\jump$?
\ed
\end{incom0}

\addtocounter{footnote}{2}\footnotetext{In a forthcoming paper, we provide the intriguing answer to 
all but the first of these questions, by showing that while there are no reals between 
$0$ and $0^\jump$, there are incomparable sets of reals between $0$ 
and $0^\jump$. This latter fact is proved by adapting the classical 
priority argument technique to the supertask context.}
 
We would like in the next few theorems to understand the relationship between
the degrees represented by a real and the degrees which are represented 
only by sets of reals.
 
\newtheorem{below}[incom]{Theorem}
\begin{below}
Any countable set of reals is infinite time computable from a real.
\label{sevenfive}
\end{below}
 
\proof  Suppose that $A=\set{z_n\mid n\in\omega}$. Let $z$ canonically 
code the sequence $\<z_n\mid n\in\omega>$. Clearly $A$ is computable from 
$z$, since an infinite time Turing machine can perform the decoding.\eth
 
\newtheorem{notbelow}[incom]{Theorem}
\begin{notbelow}
Most oracles $A\subset \R $
are not computable from a real.
\end{notbelow}
 
\proof  This is a simple counting argument. There are only $2^\omega$ 
many reals, and each of them computes at most countably many sets, so there
are only $2^\omega$ many sets which are computable from a real, but there
are $2^{2^\omega}$ many sets of reals.\eth
 
The interesting thing, however, is that one can also ensure that there are no
noncomputable reals below an oracle.
 
\newtheorem{avoid}[incom]{Theorem}
\begin{avoid}
There is a countable semi-decidable set $W$ such that any larger set 
$B\supseteq W$ computes no noncomputable reals. In fact, the
set $W$ can be taken to be the set of accidentally writable reals.
\end{avoid}
 
\proof  Let $W$ be the set of accidentally writable reals. This set is
certainly countable and semi-decidable. Moreover, it includes, and in
fact is equal to, the set of reals for which an oracle query is made
during the computation of $\varphi_p^{\R}(0)$, for any program $p$.
Notice that when $\R$ is the oracle, queries are always answered Yes.
But $W$ will also answer Yes to all those queries, since if on input
$0$ a real is written on the tape and a query made about it, it must
have been an accidentally writable real. Consequently, we have
$\varphi_p^W(0)=\varphi_p^{\R}(0)$, for any program $p$ (meaning that
if one converges, then both converge to the same answer). Furthermore,
if $W\subset B$, then similarly $\varphi_p^{B}(0)=\varphi_p^{\R}(0)$.
Now if $z$ is computable from $B$, it must be that $z=\varphi_p^{B}(0)$
for some program $p$, and consequently $z=\varphi_p^{\R}(0)$ also.
Thus, $z$ is computable from $\R$, and hence it is computable.\eth
 
\newtheorem{cone}[incom]{Corollary}
\begin{cone}
There is a set $B$ of reals which neither computes any noncomputable real, nor
is computable from any real.
\end{cone}
 
\vspace*{3in}
 
\hspace*{1in} \special{eps:Cor6-8.eps x=3in. y=3in.}
 
\proof  Simply combine the two previous arguments. There are 
$2^{2^\omega}$ many sets $B$ which contain $W$, but only $2^\omega$ many of 
them are computable from a real, so there must be some $B$ containing $W$ 
which is not computable from a real. Since it contains $W$, it also fails to
compute any noncomputable real.\eth

\newtheorem{avoid2}[incom]{Theorem}
\begin{avoid2}
For every real $x$ there is a set $A$ which computes the same reals as $x$ but
which itself is not computable from any real. 
\end{avoid2}
 
\proof  Fix the real $x$. Let $Q_x$ be the set of reals for which
a query is made during the computations of the form
$\varphi_p^{\R}(x)$.  Since the oracle is $\R$, such queries are always
answered Yes. If $B$ is a set such that $Q_x\subseteq B\subseteq \R$,
then $B$ will also give the answer Yes to all such queries, and so
$\varphi_p^{B}(x)=\varphi_p^{\R}(x)$ for any program $p$. Thus, the
reals below $B\oplus x$ are exactly those below $x$. By the counting
argument, we may find a $B$ which is not computable from any real.  Let
$A=B\oplus x$ for such a $B$. Then the reals computable from $A$ are
exactly the reals computable from $x$, and $A$ is not computable from
any real.\eth
 
\newtheorem{repr}[incom]{Theorem}
\begin{repr}
For every set $A$ there is a countable set $B\subset A$ which is 
semi-decidable in $A$ and computes the same reals as $A$.
\end{repr}
 
\proof  Let $B$ be the set of reals in $A$ for which a query is made
during a computation of the form $\varphi_p^{A}(0)$. Thus, for any
program $p$, we have $\varphi_p^{B}(0)=\varphi_p^{A}(0)$, since $B$
answers Yes whenever $A$ does for such computations. It follows that
$A$ and $B$ compute the same reals. 
The set $B$ is countable since
there are only countably many programs $p$, and each computation
$\varphi_p^{A}(0)$ mentions only countably many reals. Clearly $B$ is
semi-decidable in $A$, since one simply simulates the computations
$\varphi_p^{A}(0)$, and says Yes to any real in $A$ which appears along
the way.\eth
 
\newtheorem{high}[incom]{Jump Iteration Theorem}
\begin{high}
If a real $z$ is computable from $0^\Jump$, then so is $z^\jump$.
Indeed, we may iterate the jump operator $\jump$ many times: if $\alpha$ is
$0^\Jump$-writable, then $z^{\jump^{(\alpha)}}\leq_{\infty} 0^\Jump$.
\end{high}
 
\proof  Suppose $z=\varphi_r^{0^\Jump}(0)$. Let 
$\varphi_q(p,x)=\varphi_p^{x}(0)$. Thus, 
$p \in z^\jump\leftrightarrow\varphi_p^{z}(0)\converges\leftrightarrow\varphi_q(p,z)\converges
\leftrightarrow(q,\<p,\varphi_r^{0^\Jump}(0)>)\in 0^\Jump$. So $z^\jump$ is computable
from $0^\Jump$. Suppose now that $\alpha$ is coded by the $0^\Jump$-writable real $y$,
and we want to iterate the jump $\alpha$ many times. We use the real $y$ to 
organize the scratch tape into $\alpha$ many rows, putting the jump of each 
row into the next row (according to the order given by $y$). At limit stages, 
we write the real coding all of the earlier reals according to the organization
given by $y$. When we have finished, we have a real coding 
$z^{\jump^{(\alpha)}}$, as desired.\eth
 
\newtheorem{low}[incom]{Low Theorem}
\begin{low}
Every real below $0^\Jump$ is low. That is, if $z\ITlt0^\Jump$, then 
$z^\Jump\ITequiv 0^\Jump$. And this may be relativized to any oracle: 
every real below $A^\Jump$ is $A$-low, in the sense that if 
$z\ITlt A^\Jump$, then $z^\Jump\ITleq A^\Jump$.
\end{low}
 
\proof  Certainly $0^\Jump\ITleq z^\Jump$, so we must only show the
other direction. Suppose $z=\varphi_r^{0^\Jump}(0)$.  Let
$\varphi_q(p,x,y)=\varphi_p^{y}(x)$. Thus, 
$$(p,x)\in
z^\Jump\Iff\varphi_p^{z}(x)\converges\Iff
\varphi_q(p,x,z)\converges$$
$$\Iff(q,\<p,x,z>)\in
 0^\Jump\Iff 
(q,\<p,x,\varphi_r^{0^\Jump}(0)>)\in 0^\Jump.$$ Thus, $z^\Jump$ is
computable from $0^\Jump$ and we are done. The same argument works for
any oracle $A$.\eth
 
\newtheorem{evclosure}[incom]{Eventual Jump Theorem}
\begin{evclosure}
The class of eventually writable reals is closed under the jump
operator $\jump$. Indeed, if $z$ is an eventually writable real, and
$\alpha$ is an eventually writable ordinal, then the $\alpha^{\rm th}$
iterate $z^{\jump^{(\alpha)}}$ of the jump operator $\jump$ is still
eventually writable.
\end{evclosure}
 
\proof Suppose that $z$ is eventually written by program $p$. Consider
now the supertask algorithm that simulates $p$ on a portion of the
scratch tape. We regard the various reals that appear on the output
tape of the simulation as approximations to $z$. For every such
approximation $w$, we begin the computation to write $w^\jump$ on the
output tape, while continuing the simulation of $p$ on the scratch
tape. If we find that the approximation $w$ changes at some point, then
we disregard our previous attempts to compute $w^\jump$, and start
fresh with the new approximation. Eventually, the real $z$ will appear
as its own approximation, never subsequently to be changed, and we will
eventually write $z^\jump$ on the output tape, never subsequently to
erase it.
 
A similar argument shows that $z^{\jump^{(\alpha)}}$ is also eventually
writable. For this argument, the algorithm also computes approximations
to $\alpha$, and for each such approximation, it uses the previous
technique to compute the $\jump$ iterates of $z$. Eventually, the
algorithm will have the true approximation to $z$ and the true
approximation to $\alpha$, and after a very long time, the true
approximation to $z^{\jump^{(\alpha)}}$. So $z^{\jump^{(\alpha)}}$ is
eventually writable.\eth
 
In the next theorem, we regard a 
real as semi-decidable when there is a supertask algorithm which gives 
the affirmative answers to queries about the digits of the real.

\newtheorem{imp}[incom]{Implication Theorem}
\begin{imp}
\begin{enumerate}
\item[1.] Every writable real is semi-decidable;
\item[2.] every semi-decidable real is eventually writable;
\item[3.] every eventually writable real is accidentally writable;
\item[4.] every eventually writable real is computable from $0^\Jump$;
\item[5.] and none of these implications is reversible.
\end{enumerate}
\end{imp}
 
\proof  Clearly every writable real is semi-decidable. One simply runs
the program which writes the real, consults this real, and answers Yes 
appropriately. Also, every semi-decidable real is eventually writable. 
One simply simulates 
simultaneously for each natural number $n$ the program to semi-decide whether
$n$ is in the real. Eventually all the Yes answers will be obtained and the
real will be written on the output tape, even while it is searching for more
Yes answers, not realizing that it has already found them all. And obviously 
every eventually writable real is accidentally writable. It remains to 
check that every eventually writable real is computable from $0^\Jump$. 
Suppose $x$ is eventually written by the program $p$ (we might write in 
this case $\varphi_p(0)\uparrow=x$). Consider the algorithm which simulates
the computation of $\varphi_p(0)$. At any stage of this simulation, we can 
ask the oracle $0^\Jump$ whether the algorithm which would search for this 
snapshot to repeat again would ever find that it does. When the oracle 
answers Yes, we know we have reached the looping part of the simulation, 
and we can halt, knowing that $x$ must be written on the output tape. 
 
Now we must prove that none of the implications is reversible. The real
$0^\jump$ is semi-decidable but not writable. By the previous theorem,
the real $0^{\jump\jump}$ is eventually writable but not
semi-decidable, since it is not computable from $0^\jump$. Next, the
hard part, we must prove that there is an accidentally writable real
which is not eventually writable.  Consider the algorithm which
simulates the computation of $\varphi_p(0)$ for every program $p$, and,
at every step of these computations, writes on the output tape a real
which diagonalizes against the reals on the simulated output tapes.
Thus, for every simulated step, we write a real on the output tape
which differs from every real appearing on the simulated output tapes.
Eventually, the programs $p$ which produce eventually writable reals
have stabilized in the sense that they have reached the stage where the
real is written on their output tape, not to be subsequently changed,
and our algorithm then writes a real which is different from all of
them (and also different from some other irrelevant reals). Thus, our
algorithm writes a real which is not eventually writable. Finally, we
must show that there is a real which is computable from $0^\Jump$ which
is not eventually writable. Let us define that a computation
$\varphi_p(x)$ stabilizes if it either halts or eventually writes
a real on its output tape which is not subsequently changed. Let $S$ be
the set of programs $p$ which stabilize on input $0$. The previous
diagonal argument shows that $S$ cannot be eventually writable, since
we could have just diagonalized against the reals resulting from
programs in $S$ to arrive at the same contradiction there. But
nevertheless, the set $S$ is computable from $0^\Jump$, as we will now
show. Consider the algorithm which simulates the computation of all
$\varphi_p(0)$. Eventually, all the programs which will stabilize have
stabilized, and only then will the simulation go into an infinite
repeating loop. At any stage of the simulation, we can ask the oracle
$0^\Jump$ whether we have reached a snapshot yet which will be repeated
later. Since $0^\Jump$ can answer such questions, we will know whether
we have reached the infinite repeating loop. After having reached this
loop, we run through the loop once, checking which of the computations
$\varphi_p(0)$ change their output tape. The ones that do not are
exactly the elements of $S$. So $S$ is computable from $0^\Jump$.\eth
 
It is possible to show also that $S$ is accidentally writable. We do not 
know, however, the answer to the following question.
 
\newtheorem{acc0}[incom]{Question}
\begin{acc0}
Is every accidentally writable real computable from $0^\Jump$? Or vice 
versa? 
\end{acc0}
 
If every $0^\Jump$-clockable ordinal is also $0^\Jump$-writable, 
then we can prove that there is a real computable from $0^\Jump$ 
which is not accidentally writable.
 
\section{A new pointclass below $\underTilde\Delta^1_2$}
 
In this section we will analyze a boldface version, if you will, of the
decidable sets.  Namely, we define that a set of reals $A$ is
decidable from a real when there is a real $z$ with respect to which
$A$ is decidable.  These sets form a new natural pointclass below
$\underTilde\Delta^1_2$. The jump operator will stratify the
$\underTilde\Delta^1_2$ sets, and indeed all sets of reals, into a fine
hierarchy, ordered by $\ITleq$.
 
\newtheorem{s-alg}[incom]{Decidable Pointclass Theorem}
\begin{s-alg}
The class of sets which are decidable from a real is a $\sigma$-algebra, 
and is closed under Suslin's operation $\cal A$. Hence, they include the 
$C$-sets of Selivanovski. But this inclusion is proper, for there is a 
decidable set which is not a $C$-set.
\end{s-alg}
 
\proof  The C-sets of Selivanovski are defined to be those in the smallest $\sigma$-algebra
containing the Borel sets which is closed under Suslin's operation
$\cal A$. The sets decidable from a real are clearly closed under
complement. For the first part of the theorem, we need to prove they
are closed under countable union.  Suppose that for each natural number
$n$, the set $A_n$ is decidable by program $p_n$ from real $z_n$. Let
$z$ be a real which canonically codes $\<z_n\mid n\in \omega>$ and
$\<p_n\mid n\in\omega>$. Consider the supertask algorithm which on input
$x$ uses $z$ as an oracle to systematically check whether $x\in A_n$.
The algorithm simply simulates $p_n$ with oracle $z_n$ on input $x$,
and gets the answer yes or no whether $x\in A_n$. After our algorithm
has finished the simulations, let it output yes or no accordingly if
$x$ was in any of the sets $A_n$. Thus, $\cup_n A_n$ is decidable from
$z$, and so the class of sets decidable from reals is a
$\sigma$-algebra.
 
For the second part, recall that Suslin's operation $\cal A$ is defined
on a family of sets $\<A_s\mid s\in \omega^{<\omega}>$ to be the set of
$x$ such that there is a $y\in \omega^\omega$ such that $x\in
A_{y\restrict n}$ for every natural number $n$. Thus, $x\in{\cal
A}(\<A_s\mid s\in \omega^{<\omega}>)$ exactly when $\set{s\mid x\in
A_s}$ has an infinite branch. By intersecting the sets with the ones
preceding them, it suffices to consider only the case when $s \subset t
\rightarrow A_t \subset A_s$, so that $\set{s\mid x\in A_s}$ is a
tree.  So, suppose $A_s$ is decidable from $z_s$ by the program $p_s$
for each $s\in \omega^{<\omega}$. Let $z$ be a real which canonically
codes the other reals $\<z_s\mid s\in \omega^{<\omega}>$ as well as the
programs $\<p_s \mid s\in \omega^{<\omega}>$. Now suppose we are given
$x$, and we wish to decide if $x\in{\cal A}(\<A_s\mid s\in
\omega^{<\omega}>)$. Using $z$ as an oracle, we can simulate for each
$s$ the program $p_s$ with oracle $z_s$ on input $x$, and write out the
set $\set{s \mid x \in A_s}$. We can assume this is a tree. In $\omega$
many steps, we can write out the relation coding the Kleene-Brouwer
order on this tree, and then use the count-through algorithm to
determine if this order is a well order.  Since this occurs exactly
when the tree has no branches, this algorithm can decide whether $x\in
{\cal A}(\<A_s \mid s\in \omega^{<\omega}>)$. So the class of sets
decidable from a real is closed under Suslin's operation $\cal A$.
Thus, every $C$-set is decidable from a real.
 
Lastly, we would like to show that this inclusion is proper.  There is
a natural way to code C-sets---one just labels each node in a countable
well-founded tree with instructions for taking the union, complement,
or operation $\cal A$ of the children of that node. The leaves of the
tree are labeled with basic open sets. The set coded by such a code is
obtained by simply working up the tree, assigning a set to each node
according to the instructions on the labels. The set assigned to the
top node is the desired set. And every C-set is coded by such a code.
The labeled tree itself can easily be coded by a real, and the set of
reals coding C-set codes is $\Pi^1_1$, since the only complicated part
is that the tree must be well-founded.  Thus, the set of C-set codes is
decidable. If $w$ is a C-set code, let $A_w$ be the C-set which $w$
codes.  Let us argue that the relation $x\in A_w$ is infinite time
decidable. First, we already argued that we can decide if $w$ is in
fact a C-set code. If it is, then we can systematically decide whether
$x\in A_y$ for each code $y$ which appears as a node in the tree coded
by $w$. Thus, at the leaves of the tree, we decide if $x$ is in the
basic open set at that node. If a node instructs us to take a
complement at that node, then we flip the previous answer to its
opposite. If a node instructs us to take a union of the previous nodes,
then we search to see if $x$ is in any of the previous nodes, and write
the answer accordingly. Lastly, if a node instructs us to apply
Suslin's operation $\cal A$ to the previous nodes, then we have to
write down the Kleene-Brouwer order on the `tree' of finite sequences
coded in the children whose C-sets contain our given real $x$. If this
order is a well-order then we write No, otherwise Yes, on the node
labeled $\cal A$. At the end of this algorithm, the top node has been
labeled according to whether $x\in A_w$ or not, and we output the
answer.  It follows by this argument that the set $D=\set{w \mid w \notin
A_w}$ is decidable.  Furthermore, it is easy to see that it is not a
C-set, because it cannot be equal to $A_w$ for any code $w$. So the
theorem is proved.\eth
 
\newtheorem{lightface}[incom]{Remark}
\begin{lightface}
The previous theorem has a lightface analog. Namely, the class of 
decidable sets is effectively a $\sigma$-algebra, in the sense that
it is closed under complements and if the sets $\<A_n\mid n\in \omega>$ 
are uniformly decidable, then so is their union; and the class of 
decidable sets is closed under effective applications of $\cal A$ in the 
sense that if $\<A_s\mid s\in \omega^{<\omega}>$ is uniformly decidable, 
then so is ${\cal A}(\<A_s\mid s\in \omega^{<\omega}>)$. 
\end{lightface}
 
\vspace*{4.5in}
 
\special{eps:Thm7-2.eps x=3.6in. y=4.5in.}

The previous theorems identify a new hierarchy of pointclasses between the 
C-sets and the $\underTilde\Delta^1_2$ sets, as illustrated in the previous
diagram. One moves up the hierarchy by applying the jump operator. 
The next theorem tells us some more 
descriptive-set-theoretic information about the class of semi-decidable 
sets.  For those unfamiliar with some descriptive set theoretic terms,
\cite{Mos} provides an excellent and full description of such concepts 
as norms (p.69), the prewellordering
property (p. 200), and Spector point classes (p. 207).
 
\newtheorem{norm}[incom]{Semi-decidable Norm Theorem}
\begin{norm}
Every semi-decidable set admits a semi-decidable norm. Thus, the class
of semi-decidable sets has the pre-wellordering property, and is a 
Spector pointclass.
\end{norm} 
 
\proof  Suppose $A$ is semi-decidable by some program $p$ (for
definiteness, assume $p$ is the least such program).  Thus, $x\in A$ if
and only if $\varphi_p(x)\converges$ (and in this case the value is
$1$).  Let $\rho(x)$ be the number of steps that the computation
$\varphi_p(x)$ takes to halt, if it does halt. This is obviously a norm on
$A$. What we have to show is that the two relations $$x \leq_\rho y
\leftrightarrow x \in A \wedge [y \notin A \vee \rho(x) \leq\rho(y))]$$
and $$x <_\rho y \leftrightarrow x \in A \wedge [y \notin A \vee
    \rho(x)< \rho(y)].$$ are both semi-decidable. The first relation
$x\leq_\rho y$ holds exactly when the computation of $\varphi_p(x)$
halts before or simultaneously with $\varphi_p(y)$, and this is clearly
semi-decidable, since we can simply simulate both computations,
checking to see if $\varphi_p(x)$ halts in time.  The second relation
$x<_\rho y$ holds exactly when the computation $\varphi_p(x)$ converges
and halts strictly before $\varphi_p(y)$. This is also semi-decidable
by a similar algorithm.
 
The class of semi-decidable sets is clearly a $\Sigma$-pointclass with
the substitution property, closed under $\forall^\omega$,
$\omega$-parameterized by the programs, and, by the argument just given,
normed.  Thus, it is a Spector pointclass.\eth
 
\section{Admissibility}
 
In this section we would like to explore the connections between
infinite time computability and admissible set theory. A transitive
collection of sets is admissible when it is a model of the following
set theoretic axioms: pairing, union, infinity, Cartesian product,
$\Delta_0$-comprehension, and $\Delta_0$-collection. We will refer to
these axioms as the Kripke-Platek axioms. We refer the reader to
\cite{Bar} for an excellent account of admissibility, and content
ourselves here to say that the KP axioms form a particularly important
fragment of set theory. Though weak, this fragment is sufficiently
strong to carry out many set-theoretic constructions, such as that of
G\"odel's $L$. An ordinal $\alpha$ is admissible when $L_\alpha$ is an admissible 
set.
 
Our infinite time Turing machines of course deal only with reals, but
we can quite easily use reals to code hereditarily-countable sets.
Specifically, given a hereditarily-countable set $a$, one enumerates
$TC(\set{a})=\set{a_n\mid n\in\omega}$, and then defines a relation
$iEj\leftrightarrow a_i\in a_j$. Thus,
$\<\omega,E>\cong\<TC(\set{a}),{\in}>$, and so from $E$ we can recover
$TC(\set{a})$, and hence also $a$ (since $\{a\}$ is the unique member
of $TC(\set{a})$ which is not a member of any other member).  The
relation $E$ can be coded with a real in the usual manner. For the rest
of this section, let us assume complete familiarity with such coding
techniques, and move on now to prove the admissibility of the various
collections of sets to which the machines have access.
 
\newtheorem{admiss}{Theorem}[section]
\begin{admiss}
The class of sets coded by writable reals is admissible. 
\end{admiss}
 
\proof  The first thing to notice is that the set of reals which
are codes for sets is $\Pi^1_1$ and therefore decidable. 
To determine whether two codes actually code the same set
is $\Sigma^1_1$ (since this is true if and only if there is an
isomorphism between the indexed sets) and therefore
also decidable, as is determining whether codes $x$ and $y$ 
code sets such that $A_x\in A_y$.
Thus, it is easy to see that $\Delta_0$ facts about the sets coded by
reals are decidable from the codes. Furthermore, from a code $z$ for
the set $A_z$, a code can be generated for any element $a\in A_z$,
knowing merely which natural number $n$ represents $a$ in the code
$z$.  The class of sets coded by writable reals is clearly closed under
union, pairing, and difference. And it satisfies
$\Delta_0$-comprehension because from a code $z$ for a set $A_z$, and a
$\Delta_0$-formula $\psi$, we can systematically compute a code for the
set $\set{a\in A_z\mid \psi(a)}$ by simply computing whether the formula
$\psi$ holds separately for each member of the set, and then putting
all desired elements together into a code. Finally, we will prove that
$\Delta_0$-collection holds. So suppose $z$, coding the set $A_z$, is
writable, and for every $a\in A_z$ there is a set $b$ coded by a
writable real such that $\psi(a,b)$ holds, where $\psi$ is a $\Delta_0$
formula.  What we need to find is a set $B$ which is coded by a
writable real such that for every $a\in A_z$ there is a witness $b\in
B$ such that $\psi(a,b)$.  Let us now show there is such a $B$. Fix
$z$, and consider the supertask algorithm which first writes $z$ on a portion of
the scratch tape. Now, we will slowly write down the code for $B$, by
searching for witnesses. Every element of $A_z$ is indexed by some
natural number in the code $z$, and for each such element we will start
the algorithm which simulates the computation of $\varphi_p(0)$ for
every program $p$, until a witness is produced for the given element.
By hypothesis, we will eventually find a $y=\varphi_p(0)$ for some
program $p$ such that $\psi(a,A_y)$ holds, and we copy the code $y$ to
represent an element of the set $B$ we are building.  After doing this
for each element of the set coded by $z$, we have written a code for
the set $B$, and we may halt. So $\Delta_0$-collection holds, and
therefore the set of writable reals is a model of KP, and hence
admissible.\eth
 
\newtheorem{admiss2}[admiss]{Corollary}
\begin{admiss2}
The supremum $\lambda$ of the writable ordinals is an admissible ordinal.
\end{admiss2}
 
\proof  Every writable real is constructible because we could just run the
    computation in $L$. Also, if $\alpha$ is a writable ordinal,
    then we claim that $L_\alpha$ is coded by a writable real.
This is true because, as in the Lost Melody Theorem \ref{melody}, an infinite
time Turing machine can simulate the
    construction of $L_\alpha$ given a code for $\alpha$; given a
    code for $L_\beta$, one obtains a code for $L_{\beta+1}$ by
    enumerating the definitions and systematically writing down the
    codes for the definable subsets of $L_{\beta}$; and at limit stages
    of the $L$ construction, given the sequence of codes for the
    earlier stages, one simply writes down the code for the union.
    Thus, if $\lambda$ is the supremum of the writable ordinals,
    $L_\lambda$ is a subclass of the collection of sets coded by
    writable reals (it is not clear whether these collections are
    distinct). Furthermore, the proof of the previous theorem can be
    modified to use $L_\lambda$ rather than the class of all sets
    coded by writable reals. Thus, rather than just searching for
    writable reals coding hereditarily countable sets, one searches
    for writable reals coding writable ordinals $\alpha$, and then
    uses these ordinals to construct $L_\alpha$ in the manner we
    have just explained, and then searches for the witnesses in
    these $L_\alpha$. Thus, $L_\lambda$ is admissible, and
    consequently $\lambda$ is admissible. \eth
 
The next theorem will show that $\lambda$ is quite high up in the
hierarchy of admissible ordinals. An ordinal $\beta$ is
{\it recursively inaccessible} when it is an admissible limit of admissible
ordinals. The ordinal $\beta$ is {\it indescribable} by a class of properties
if there is a real $x$ coding $L_\beta$ such that for any property in the 
class which is true of $x$ there is an $\alpha<\beta$ and a real
$y$ coding $L_\alpha$ having the very same property. 
 
\newtheorem{admiss3}[admiss]{Indescribability Theorem}
\begin{admiss3}
The supremum $\lambda$ of the writable ordinals is recursively
inaccessible.  Indeed, it is the $\lambda^{\rm th}$ recursively
inaccessible ordinal, and the $\lambda^{\rm th}$ such fixed point,
and so on. This is because $\lambda$ is indescribable by 
$\Pi^1_1$ properties. Indeed, $\lambda$ is indescribable by 
semi-decidable properties.
\label{indescribe}
\end{admiss3}
 
\proof  Let $\lambda$ be the supremum of the writable ordinals. We have
already shown that $\lambda$ is admissible. If it is not a limit of
admissible ordinals, then there is some largest admissible
$\delta<\lambda$, which is consequently writable. Now, whether a real
codes an admissible ordinal is infinite time decidable, since from the
code for an ordinal $\alpha$ we have already explained how to get a
code for the set $L_\alpha$, and it is decidable to check whether the
set coded by a given real is admissible, since this consists in merely
checking that the real codes a model of a certain recursively
axiomatized theory. So, consider the supertask algorithm which first writes a
real coding $\delta$ on a portion of the scratch tape, and then
simultaneously simulates the computation of $\varphi_p(0)$ for every
program $p$. For every real appearing during these computations, the
algorithm checks to see if it codes an ordinal larger than $\delta$
which is admissible. If so, the algorithm gives that real as the output
and halts.  Since $\lambda$ is accidentally writable, this algorithm is
bound to find, and write, the code for an admissible ordinal above
$\delta$. By our assumption, this ordinal must be at least $\lambda$,
contradicting the fact that $\lambda$ is the supremum of the writable
ordinals. Thus, $\lambda$ must be an admissible limit of admissible
ordinals.
 
By a similar argument we will now show that it cannot be the least such
limit.  Suppose $\lambda$ is the least limit of admissible ordinals above
$\delta$. Now, it is a decidable question whether a real codes an
ordinal which is an admissible limit of admissible ordinals; one simply
tests first whether the real codes an admissible ordinal, and then
tests whether there is any index for a smaller ordinal with no
admissible ordinals in-between. So consider the supertask algorithm which first
writes a code for $\delta$ on a portion of the scratch tape, and then
searches for an accidentally writable real which codes an admissible
limit of admissible ordinals, and tests if it is larger than $\delta$.
When such a real is found, the algorithm gives it as the output and
halts. By hypothesis, this algorithm will write a real at least as
large as $\lambda$, a contradiction.
 
Suppose now that $\lambda$ is the $\delta^{th}$ admissible ordinal for
some $\delta<\lambda$. Consider the algorithm which first writes
$\delta$ on a portion of the scratch tape, and then searches for
admissible ordinals which are coded by accidentally writable reals.
Each time one is found which is larger than the previous ones, it is
written on a portion of the scratch tape labeled with the index of the
next element in the relation coding $\delta$.  When every element of
the field of that relation is taken care of, the algorithm writes the
corresponding ordinal on the tape, and halts. By assumption, this real
codes an ordinal at least as large as $\lambda$, a contradiction.
 
The indescribability argument is no different. Let $x$ be an eventually
writable real which codes $L_\lambda$ in the manner of which we have
been speaking. If $x$ has some semi-decidable property, then there must be a
smaller ordinal $\alpha<\lambda$ with an accidentally writable real $y$
coding $L_\alpha$ with the very same property, since otherwise the
algorithm which went searching for such a real would be able to write
an ordinal larger than $\lambda$, which is impossible. In fact, there
must be a writable $y$ with that property, since the algorithm will halt
when one is found. Consequently, $\lambda$ is
indescribable by semi-decidable properties.\eth
 
\newtheorem{Qleast}[admiss]{Question}
\begin{Qleast}
Is $\lambda$ the least ordinal which is indescribable by semi-decidable properties?
\end{Qleast}
 
Next, we will prove the corresponding facts about the class of sets
coded by eventually writable reals.
 
\newtheorem{admiss4}[admiss]{Theorem}
\begin{admiss4}
The class of sets coded by eventually writable reals is admissible.
\end{admiss4}
 
\proof  This proof is very similar to the corresponding proof for
writable reals. But it also has the flavor of a finite injury priority 
argument in classical recursion theory. The class in question 
is closed under the rudimentary 
functions, and satisfies $\Delta_0$-comprehension just as before. The 
hard part is $\Delta_0$-collection. So suppose $z$ is eventually writable,
and for every $a\in A_z$ there is a set $b$ coded by an eventually writable 
real such that $\psi(a,b)$, where $\psi$ is some fixed $\Delta_0$ formula. 
We have to collect witnesses into a set which is eventually writable. Consider
the supertask 
algorithm which writes approximations to $z$ on a portion of the scratch 
tape. We will use these approximations to $z$ to write down a code for a
set of witnesses. Eventually the correct approximation will be written, and
we will eventually write down a code for a set of witnesses. So, for 
each approximation $w$ to $z$, we start writing down the code for a set of
witnesses that will work for $w$, simultaneously computing better 
approximations to $z$. For each index $n$, representing a set $a$ in $A_w$, we 
search for an eventually writable witness $b$ for $a$. We do this by 
simulating for every program $p$ the computation of $\varphi_p(0)$ and 
testing the approximate outputs $b$ for these computations to see if 
$\psi(a,b)$ holds. If so, we copy the code for $b$ to be an element of $B$
on the true output tape. Periodically, however, the approximation $w$ may 
change because it has not yet stabilized, and when this occurs, we erase 
the witness $b$ from the output tape, and move to the 
approximation produced by another program, in such a way that we give every
program a chance to produce the correct witness. Eventually, 
by hypothesis, we will hit on a witness which has stabilized. Furthermore, 
we do this simultaneously for each set $a\in A_w$, and every time the 
approximation $w$ changes, we start the process completely over with the 
new approximation. Eventually, $z$ will be written as its own approximation, 
and the witness computations will hit on the stabilizing witness which 
will be copied as elements of the set $B$ on the output tape. While the
algorithm will continue searching for better approximations, it will never
change the set $B$ after this point, and so we may collect witnesses into
an eventually writable set $B$. So $\Delta_0$-collection holds.\eth
 
Let us define that a set $A$ is {\it eventually} decidable when there is a 
program $p$ such that for any $x$ the computation $\varphi_p(x)$ never halts, 
but eventually has either a $1$ or a $0$ written on the output tape (never 
subsequently to change), respectively, depending on whether $x\in A$ or not. 
It is clear that every semi-decidable set is eventually decidable.
Similarly, $A$ is eventually semi-decidable when there is a program which
for any $x$ eventually writes the Yes answers on the output tape, accordingly, 
depending on whether $x\in A$ or not. 
 
\newtheorem{admiss5}[admiss]{Corollary}
\begin{admiss5}
The supremum $\zeta$ of the eventually writable ordinals is an
admissible ordinal, a recursively inaccessible ordinal, the $\zeta^{\rm
th}$ recursively inaccessible ordinal, the $\zeta^{\rm th}$ such fixed
point, and so on. It is indescribable by decidable properties and even
by eventually semi-decidable properties.
\end{admiss5}
 
\proof  Let $\zeta$ be the supremum of the eventually writable ordinals.
    By the argument of \ref{indescribe} it follows that $L_\zeta$ is a subclass
    of the sets coded by eventually writable reals (though again it
    is not clear whether these classes are distinct). In any case,
    the proof of the previous theorem is easily adapted as in \ref{indescribe},
    by limiting the searches to $L_\zeta$ rather than any eventually
    writable set, to show that $L_\zeta$ is admissible.
    Consequently, $\zeta$ is an admissible ordinal. 

Next, let us argue that $\zeta$ is accidentally writable, even though
it is not eventually writable. Consider the supertask algorithm which simulates
the computation of $\varphi_p(0)$ simultaneously for every program $p$.
Our algorithm checks which of the output approximations code well
orders, and writes a code for the sum of the ordinals coded by those
that do. Eventually, the eventually writable ordinals are written as
their own approximations, and so we have written an ordinal bigger than
every eventually writable ordinal. Since the accidentally writable
ordinals are closed under initial segment, it follows that $\zeta$
itself is accidentally writable.
 
Now we can show that $\zeta$ is recursively inaccessible, and
indescribable, and so on, just as before. There is an accidentally
writable real $x$ coding $L_\zeta$. Suppose that $x$ has some
eventually semi-decidable property.  Consider the algorithm which
searches for reals with this very same property, and writing the real
coding the corresponding ordinal.  If $\zeta$ is the least ordinal such
that $L_\zeta$ is coded by an accidentally writable real with that
property, then this algorithm will eventually settle on a real coding
an ordinal at least as large as $\zeta$, contradicting the fact that
$\zeta$ is the supremum of all eventually writable ordinals. Thus,
there must be a smaller $\alpha<\zeta$ with an accidentally writable real
$y$ coding $L_\alpha$. The real $y$ will be eventually writable, since
the algorithm will settle on $y$. So $\zeta$ is indescribable by
eventually semi-decidable properties.\eth
 
In our penultimate theorem, let us give necessary and sufficient
conditions on the question of whether every clockable ordinal is
writable. We have gone back and forth on this issue, and we simply do
not know the answer. As usual, $\lambda$ is the supremum of the
writable ordinals, and $\gamma$ is the supremum of the clockable
ordinals.
 
\newtheorem{neccsuff}[admiss]{Theorem}
\begin{neccsuff}
The following are equivalent:
\begin{enumerate}
\item[1.] $\lambda=\gamma$.
\item[2.] Every clockable ordinal is writable.
\item[3.] $\lambda$ is a limit of clockable ordinals.
\item[4.] The halting problem $h_\lambda$ is not decidable.
\item[5.] $\gamma$ is admissible.
\end{enumerate}
\end{neccsuff}
 
\proof $(1\leftrightarrow 2)$ This is clear, since there are no gaps in the
writable ordinals, and certainly $\lambda\leq\gamma$, since $\lambda$
is the order-type of the clockable ordinals.
 
$(1\leftrightarrow 3)$ The forward direction is clear, since $\gamma$
is a limit of clockable ordinals. For the converse direction, suppose
that $\lambda$ is a limit of clockable ordinals, but $\lambda<\gamma$.
Thus, $\lambda$ begins a gap in the clockable ordinals
(recall that $\lambda$ itself is not clockable by \ref{threetwelve}). 
Let $\bar\lambda$ be the supremum of the ordinals which are writable in
less than $\lambda$ many steps of computation. The lengths of these
computations must in fact be unbounded in $\lambda$, since otherwise
they would be bounded by some clockable ordinal $\beta<\lambda$, and we
could run them all at once, with a $\beta$-clock, and write a real
coding an ordinal larger than $\bar\lambda$ in fewer than $\lambda$
many steps, contradicting the definition of $\bar\lambda$. Also, by
essentially the same argument, it must be that $\bar\lambda<\lambda$
since otherwise we could run all the computations up to some clockable
$\beta$ beyond $\lambda$ and write a real coding an ordinal larger than
$\lambda$, a contradiction. Now we can complete the argument by
observing that the map which takes $\alpha<\bar\lambda$ to the length
of the shortest computation which writes $\alpha$ is an unbounded map
from $\bar\lambda$ to $\lambda$. Since it is $\Sigma_1$-definable in
$L_\lambda$, where the computations exist, this contradicts the fact
that $\lambda$ is admissible.
 
$(1\leftrightarrow 4)$ The forward direction again is clear since $h$ is not
decidable.  For the converse, observe that if $\lambda$ is below some
least clockable ordinal $\beta$, then $h_\lambda=h_\beta$, and
$h_\beta$ is decidable by Theorem \ref{seven}.
 
$(1\leftrightarrow 5)$ The forward direction follows since $\lambda$ is
admissible.  For the converse, suppose that $\lambda<\gamma$. By the
argument showing $(3\rightarrow 1)$, it must be that the lengths of the
computations which write the various ordinals below $\lambda$ are
unbounded in $\gamma$. Consequently, in $L_\gamma$ we have a
$\Sigma_1$-definable map from $\lambda$ unbounded in $\gamma$. So
$\gamma$ cannot be admissible.\eth
 
The argument in $(1\leftrightarrow 3)$ shows that if $\lambda<\gamma$, then there 
are no ordinals in the interval $[\lambda,\gamma]$ which are both 
admissible and a limit of clockable ordinals. 
 
We would like to conclude our paper by answering a question we teased the 
reader with way back in section three. Namely, is the ordinal $\wCK$ 
clockable? The answer, by the following theorem, is No.
 
\newtheorem{noadmiss}[admiss]{Theorem}
\begin{noadmiss}
No admissible ordinal is clockable.
\end{noadmiss}
 
\proof  Suppose that $\alpha$ is a clockable limit ordinal. We will
show that $\alpha$ is not admissible by showing that there is a
function $f:\omega\to\alpha$, unbounded in $\alpha$, which is
$\Sigma_1$-definable in $L_\alpha$. It follows that $L_\alpha$ cannot
model $\Sigma_1$-collection, and therefore cannot be admissible. It
suffices to consider the case when $\alpha$ is a limit of limit
ordinals. Let $p$ be a program which on input $0$ halts in exactly
$\alpha$ steps. And let us pay attention to exactly the manner in which
the program halts at stage $\alpha$. At stage $\alpha$ the head is at
the extreme left of the tapes and in the {\it limit} state.
Furthermore, the values appearing in the three cells under the head
have caused the machine to halt. Thus, this must be the first time at
a limit ordinal that those cells have that particular pattern, or the
machine would have halted earlier. The 0s on the tape
must have been $0$ from some point on, and perhaps
one of the $1$s was also like that. But there must have been a $1$ on
the cell under the head which had alternated infinitely often before
$\alpha$, since otherwise the snapshot at $\alpha$ would have been
obtained at one of the limit stages before $\alpha$ (since by
assumption such stages are unbounded in $\alpha$), and the program
would have halted earlier. Furthermore, stage $\alpha$ must be the
first stage by which that cell had alternated from $1$ to $0$ and back
infinitely many times, since otherwise again the program would have
halted earlier. Let $h(n)$ be the stage at which the cell alternated
the $n^{\rm th}$ time. We have argued that these stages are unbounded
in $\alpha$. Furthermore, the function $h$ is $\Sigma_1$-definable in
$L_\alpha$, since the initial segments of the computation of
$\varphi_p(0)$ all live in $L_\alpha$, and are defined by a $\Sigma_1$
definition there.  Thus, over $L_\alpha$ there is a $\Sigma_1$
definable map from $\omega$ unbounded in $\alpha$, and so $\alpha$ is
not admissible.\eth

\noindent The authors can be reached at the following addresses:

\vspace*{.1in}

\noindent Joel David Hamkins, Mathematics, City University of New York, 
\newline \hspace*{.35in} College 
of Staten Island, Staten Island, NY 10314; 
\newline \hspace*{.35in} hamkins@postbox.csi.cuny.edu

\vspace*{.1in}

\noindent Andy Lewis, Mathematics, Virginia Commonwealth University,
Box \#842014, \newline \hspace*{.35in} Richmond, Va. 23284-2019; amlewis@saturn.vcu.edu

\end{document}